\documentclass[a4paper,12pt]{amsart}

\newcommand{\ind}{\text {\rm ind}}
\newcommand{\res}{\text {\rm res}}
\newcommand{\supp}{\text {\rm supp}}

\newcommand{\Ad}{\text {\rm Ad}}
\newcommand{\car}{\mathcal R}
\newcommand{\QQ}{\mathbf Q}
\newcommand{\bbq}{\bar {\QQ }}

\def\ge{\geqslant}
\def\le{\leqslant}
\def\a{\alpha}

\def\g{\gamma}

\def\d{\delta}
\def\D{\Delta}

\def\e{\epsilon}

\def\o{\omega}

\def\s{\sigma}

\def\th{\theta}

\def\i{^{-1}}

\def\cc{\mathcal C}
\def\cd{\mathcal D}
\def\ce{\mathcal E}
\def\cf{\mathcal F}

\def\co{\mathcal O}
\def\cp{\mathcal P}

\def\car{\mathcal R}
\def\cs{\mathcal S}

\def\cw{\mathcal W}
\def\cz{\mathcal Z}

\def\tz{\tilde Z}

\def\tx{\tilde X}

\def\ss{\text ss}

\usepackage{amssymb} 
\usepackage{latexsym} 
\usepackage{amsfonts} 
\usepackage{amsmath} 
\usepackage{eucal} 
\usepackage{bm} 
\usepackage{bbm} 
\usepackage{graphicx} 
\usepackage[english]{varioref} 
\usepackage[nice]{nicefrac} 
\usepackage[all]{xy}

\newcommand{\kk}{\Bbbk}

\newcommand{\sgn}{\textrm{sgn}}

\usepackage{amsthm}

\theoremstyle{plain}
\newtheorem{thm}{Theorem}[section] 
\newtheorem*{thm*}{Theorem} 
 \newtheorem{prop}[thm]{Proposition}
 \newtheorem{lem}[thm]{Lemma}
 \newtheorem{cor}[thm]{Corollary}
 
 \newtheorem{conj}[thm]{Conjecture}

\theoremstyle{definition}

\theoremstyle{remark}

\newtheorem*{remark*}{Remark}
\newtheorem*{claim*}{Claim}

\begin{document}

\author{Xuhua He}
\address{Department of Mathematics, The Hong Kong University of Science and Technology,
Clear Water Bay, Kowloon, Hong Kong}
\email{maxhhe@ust.hk}
\thanks{The author is partially supported by partially supported by HKRGC grants 601409 and DAG08/09.SC03.}
\title[]{Character sheaves on the semi-stable locus of a group compactification}

\begin{abstract}
We study the intermediate extension of the character sheaves on an adjoint group to the semi-stable locus of its wonderful compactification. We show that the intermediate extension can be described by a direct image construction. As a consequence, we show that the ``ordinary'' restriction of a  character sheaf on the compactification to a boundary piece inside the semi-stable locus is a shift of semisimple perverse sheaf and is closely related to Lusztig's restriction functor (from a character sheaf on a reductive group to a direct sum of character sheaves on a Levi subgroup). We also provide a (conjectural) formula for the boundary values inside the semi-stable locus of an irreducible character of a finite group of Lie type, which gives a partial answer to a question of Springer \cite{Sp2}. This formula holds for Steinberg character and characters coming from generic character sheaves. In the end, we verify Lusztig's conjecture \cite[12.6]{L3} inside the semi-stable locus of the wonderful compactification. 
\end{abstract}

\maketitle

\section*{Introduction}

\subsection{} Let $G$ be a connected, semisimple algebraic group of adjoint type over an algebraically closed field $k$. In \cite{L3}, Lusztig introduced a decomposition of the wonderful compactification $\bar{G}$ of $G$ into $G$-stable pieces. The group $G$ itself is a $G$-stable piece and each $G$-stable piece is a smooth, locally closed subvariety of $\bar{G}$ and the $G$-orbits on each piece (for the diagonal $G$-action) naturally correspond to the ``twisted'' conjugacy classes of a smaller group. Moreover, this correspondence leads to a natural equivalence between the bounded derived category of $G$-equivariant, constructible sheaves on that piece and the bounded derived category of certain constructible sheaves on the smaller group that are equivariant under the ``twisted'' conjugation action (see \cite[12.3]{L3}). 

Character sheaves on a reductive group are some special simple perverse sheaves on the group that are equivariant under the (``twisted'') conjugation action. The theory of character sheaves was developed by Lusztig in the series of papers \cite{L1} (for conjugation action) and \cite{L2} (for ``twisted'' conjugation action). Now using the natural equivalence we discussed above, one can define the character sheaves on each $G$-stable piece. The character sheaves on $\bar{G}$ are the intermediate extensions to $\bar{G}$ of the character sheaves on the $G$-stable pieces (see \cite[12.3]{L3}). The most interesting cases are the intermediate extension to $\bar{G}$ of the character sheaves on $G$. Roughly speaking, these sheaves can be regarded as the objects that describe the behavior at infinity of the character sheaves on $G$. 

\subsection{} In order to understand the intermediate extensions to $\bar{G}$ of the character sheaves on a $G$-stable piece, in \cite{H2} we gave a second definition of character sheaves on $\bar{G}$ by imitating the definition of character sheaves on groups. This new definition coincides with Lusztig's definition we mentioned in the previous subsection (see \cite[Corollary 4.6]{H2}). Moreover, using the new definition, one can show that the character sheaves on $\bar{G}$ have the following nice property (see \cite[Section 4]{H2}):

Let $i$ be the inclusion of a $G$-stable piece to $\bar{G}$, then

(1) for any character sheaf $C$ on $\bar{G}$, any perverse constituent of $i^* (C)$ is a character sheaf on that piece;

(2) for any character sheaf $C$ on that piece, any perverse constituent of $i_! (C)$ is a character sheaf on $\bar{G}$. 

\subsection{} However, analyzing the intermediate extension of a character sheaf on a $G$-stable piece is still a challenging problem. In \cite{Sp2}, Springer listed some interesting questions in this direction. One interesting question is to study the boundary values of an irreducible character of a finite group of Lie type.

A technical difficulty in analyzing the intermediate extension is as follows.

A character sheaf on $G$ can be understood in terms of ``admissible complex'', which is obtained by pushing forward of some intersection cohomology complex under some small, proper map to the closure of a Lusztig's stratum of $G$. 

Using the $G$-stable piece decomposition of $\bar{G}$ and the natural correspondence between the $G$-stable pieces and the smaller groups, one is able to generalize Lusztig's stratification on $G$ to a decomposition on $\bar{G}$. However, an explicit description of the closure to $\bar{G}$ of  a Lusztig's stratum is still unknown. A more serious problem is that the small map we used to construct ``admissible complex'' on $G$ doesn't extend to a small map on $\bar{G}$.  

\subsection{} In this paper, we will study the intermediate extension of a character sheaf on $G$, not to $\bar{G}$, but to the semi-stable locus $\bar{G}^{ss}$ of $\bar{G}$, an open smooth subvariety of $\bar{G}$ that contains $G$. In fact, $\bar{G}^{ss}$ is a union of some $G$-stable pieces. An explicit description of $\bar{G}^{ss}$ was obtained in an joint work with Starr \cite{HS}. 

The idea of studying intermediate extension to $\bar{G}^{ss}$ instead of $\bar{G}$ comes from geometric invariant theory. Now we make a short digression from character sheaves and discuss about some basic ideas in the theory of geometric invariant theory. 

Let $H$ be a linear algebraic group and $X$ be a $H$-variety. When considering the quotient space, a main problem is that the quotient $X/H$ may not exists in the category of algebraic varieties. Geometric invariant theory suggests a method to distinguish ``good'' $H$-orbits from ``bad'' $H$-orbits in the sense that the union of ``good'' $H$-orbits form an open subvariety $U$ of $X$ and $U/H$ exists. 

Motivated by this, one may wonder if the ``good'' $G$-orbits on $\bar{G}$ are still good in the study of character sheaves in the sense that the intermediate extension of a character sheaf on $G$ to the union of ``good'' orbits can be analyzed. The answer is YES and this is what we are going to do in this paper. 

\subsection{} Now let us consider the closure of a Lusztig's stratum in $\bar{G}$. If we take the limit in the direction of unipotent elements in the stratum, then by the results in \cite{H1} and \cite{HT1}, the boundary points are outside the semi-stable locus. On the other hand, taking the limit in the direction of semisimple elements in the stratum is more or less the same as calculating the closure of some subvariety in a toric variety. This naive thought suggests that the closure to $\bar{G}^{ss}$ of a Lusztig's stratum can be described explicitly. 

The explicit description will be obtained in section 3. Moreover,  the small map we used to construct ``admissible complex'' on $G$ extends to a small map on $\bar{G}^{ss}$. Based on this result, the intermediate extension of an ``admissible complex'' to $\bar{G}^{ss}$ can also be described by a direct image construction. This is a generalization of \cite[Proposition 5.7]{L2}. 

Moreover, the restriction of the direct image to a boundary piece inside the semi-stable locus can be calculated explicitly and is closely related to Lusztig's restriction functor introduced in \cite[3.8]{L1} and \cite[23.3]{L2}. The precise statement can be found in Theorem \ref{ppp}. Based on this, we give a (conjectural) formula for the boundary values inside the semi-stable locus of a character of a finite group of Lie type. The formula is true if the (virtual) character is obtained from the direct image construction. This gives a partial answer to a question of Springer \cite[Problem 10]{Sp2}. 

\subsection{} There is a special character sheaf $S$ on $G$ that characterizes the semisimple elements of $G$. This sheaf is the alternating sum of the induced sheaves from the trivial local systems on the standard parabolic subgroups of $G$. In \cite[12.6]{L3}, Lusztig generalized the notion of semisimple elements to $\bar{G}$ and conjectured that the intermediate extension to $\bar{G}$ of this sheaf characterizes the semisimple elements of $\bar{G}$. 

It is known that the semisimple elements of $\bar{G}$ lie in the semi-stable locus. We will calculate the intermediate extension of $S$ to $\bar{G}^{ss}$ and verify Lusztig's conjecture inside the semi-stable locus. 

In order to do this, we will consider the intermediate extension of the induced sheaf from the trivial local system on a standard parabolic subgroup $P$. Therefore we need to understand the closure of $P$ in $\bar{G}^{ss}$ and the intermediate extension of trivial local system on $P$ to this closure. 

Let $B$ be a Borel subgroup of $P$. Then $P$ is stable under the action of $B \times B$ and the closure of $P$ in $\bar{G}$ was obtained in \cite[Corollary 2.5]{Sp1} in terms of the union of certain $B \times B$-orbits. However, $\bar{G}^{ss}$ is not stable under the action of $B \times B$. To describe the closure of $P$ in $\bar{G}^{ss}$, we have to use the $P$-stable pieces, introduced by Lu and Yakimov as a generalization of the notation of $B \times B$-orbits and $G$-stable pieces.  Although the closure of $P$ in $\bar{G}$ is not smooth in general, the closure of $P$ in $\bar{G}^{ss}$ is always smooth. Therefore, the intermediate extension of trivial local system on $P$ to the closure of $P$ in $\bar{G}^{ss}$ is just the trivial local system on that closure. Now we can explicitly calculate the intermediate extension of $S$ to $\bar{G}^{ss}$. 

\subsection{} We now review the content of this paper in more detail.

In section 1, we recall the definition and properties of $P$-stable pieces. In section 2, we give an explicit description of  the closure of a parabolic subgroup in $\bar{G}^{ss}$ and prove that the closure is smooth. In section 3, we obtain the closure of a Lusztig's stratum of $G$ in $\bar{G}^{ss}$. In section 4, we study the intermediate extension of  a character sheaf on $G$ to $\bar{G}^{ss}$ and verify Lusztig's conjecture inside $\bar{G}^{ss}$.

\section{$\car$-stable pieces on the wonderful compactification}

\subsection{} Let $G$ be a connected reductive algebraic group over an algebraically closed field $k$. 
Let $B$ be a Borel subgroup of $G$, $T \subset B$ be a maximal torus and $B^-$ be the opposite Borel subgroup. Let $I$ be the set of simple roots and $W=N_G(T)/T$ be the corresponding Weyl group. For any $w \in W$, we choose a representative $\dot w$ of $w$ in $N_G(T)$. 

For $J \subset I$, let $W_J$ be the subgroup of $W$ corresponding to $J$ and $W^J$ (resp. $^J W$) be the set of minimal length coset representatives of $W/W_J$ (resp. $W_J \backslash W$). Let $w^J_0$ be the unique element of maximal length in $W_J$. (We simply write $w_0$ for $w^I_0$.) For $J, K \subset I$, we write $^J W^K$ for $^J W \cap W^K$.

For $J \subset I$, let $\Phi_J$ be the set of roots that are linear combination of simple roots in $J$. Let $P_J \supset B$ be the standard parabolic subgroup defined by $J$ and $P^-_J \supset B^-$ be the opposite of $P_J$. Let $L_J=P_J \cap P^-_J$ and $G_J=L_J/Z(L_J)$. For any parabolic subgroup $P$, we denote by $U_P$ its unipotent radical and $H_P$ the inverse image of the connected center of $P/U_P$ under $P \to P/U_P$. We simply write $U$ for $U_B$ and $U^-$ for $U_{B^-}$.

For any $g \in G$ and subvariety $H \subset G$, we write ${}^g H$ for $g H g \i$. 

\

Now we will review the $\car$-stable pieces introduced in \cite{LY}. We will follow the approach in \cite{H4}. 

\subsection{} A triple $c=(J_1, J_2, \d)$ consisting of $J_1, J_2 \subset I$ and an isomorphism $\d: W_{J_1} \to W_{J_2}$ with $\d(J_1)=J_2$ is called an {\it admissible triple} of $W \times W$. For an admissible triple $c=(J_1, J_2, \d)$, set $W_c=\{(w, \d(w)); w \in W_{J_1}\} \subset W \times W$. 

Let $c=(J_1, J_2, \d)$ and $c'=(J'_1, J'_2, \d')$ be admissible triples. For $w_1 \in W^{J_1}$ and $w_2 \in {}^{J'_2} W$, set \begin{gather*} I(w_1, w_2, c, c')=\max\{K \subset J_1; w_1(K) \subset J'_1 \text{ and } \d' w_1(K)=w_2 \d(K)\}, \\ [w_1, w_2, c, c']=W_{c'} (w_1 W_{I(w_1, w_2, c, c')}, w_2) W_c \subset W \times W. \end{gather*}

Then $W \times W=\sqcup_{w_1 \in W^{J_1}, w_2 \in {}^{J'_2} W} [w_1, w_2, c, c']$. See \cite[Proposition 2.4 (1)]{H4}. 

Moreover, define an automorphism $\s: W_{I(w_1, w_2, c, c')} \rightarrow W_{I(w_1, w_2, c, c')}$ by $\s(w)=\d \i \bigl(w_2 \i \d'(w_1 w w_1 \i) w_2 \bigr)$. Then map $W_{I(w_1, w_2, c, c')} \rightarrow W_1 \times W_2$ defined by $w \rightarrow (w_1 w, w_2)$ induces a bijection from the $\s$-twisted conjugacy classes on $W_{I(w_1, w_2, c, c')}$ to the double cosets $W_{c'} \backslash [w_1, w_2, c, c'] /W_c$. See \cite[Proposition 2.4 (2)]{H4}. 

Let $\co$ be a double coset in $W_{c'} \backslash (W \times W) /W_c$. Then $\co \cap (W^{J_1} \times {}^{J'_2} W)$ contains at most one element (see \cite[Corollary 2.5]{H4}). If $\co \cap (W^{J_1} \times {}^{J'_2} W) \neq \emptyset$, then we call $\co$ a distinguished double coset. We denote by $\co_{\min}$ the set of minimal length elements in $\co$. We have a natural partial order on the set of distinguished double cosets defined as follows: $\co \le \co'$ if for some (or equivalently, any) $w' \in \co'_{\min}$, there exists $w \in \co_{\min}$ with $w \le w'$. See \cite[4.7]{H4}. 

\subsection{}\label{cc} An {\it admissible triple} of $G \times G$ is by definition a triple $\cc=(J_1, J_2, \th_{\d})$ consisting of $J_1, J_2 \subset I$, an isomorphism $\d: W_{J_1} \to W_{J_2}$ with $\d(J_1)=J_2$ and an isomorphism $\th_{\d}: L_{J_1} \to L_{J_2}$ that maps $T$ to $T$ and the root subgroup $U_{\a_i}$ (for $i \in J_1$) to the root subgroup $U_{\a_{\d(i)}}$. Then an admissible triple $\cc=(J_1, J_2, \th_{\d})$ of $G \times G$ determines an admissible triple $c=(J_1, J_2, \d)$ of $W \times W$. For an admissible triple $\cc=(J_1, J_2, \th_{\d})$, define $$\car_{\cc}=\{(p, q); p \in P_{J_1}, q \in P_{J_2}, \th_{\d} (\bar p)=\bar q\},$$ where $\bar p$ is the image of $p$ under the map $P_{J_1} \to L_{J_1}$ and $\bar q$ is the image of $q$ under the map $P_{J_2} \to L_{J_2}$. 

Let $\cc=(J_1, J_2, \th_{\d})$ and $\cc'=(J'_1, J'_2, \th_{\d'})$ be admissible triples. For $w_1 \in W^{J_1}$ and $w_2 \in {}^{J'_2} W$, set $$[w_1, w_2, \cc, \cc']=\car_{\cc'} (B \dot w_1 B, B \dot w_2 B) \car_{\cc} \subset G \times G.$$

For any distinguished double coset $\co \in W_{c'} \backslash (W \times W)/W_c$, we also write $[\co, \cc, \cc']$ for $[w_1, w_2, \cc, \cc']$, where $(w_1, w_2)$ is the unique element in $\co \cap (W^{J_1} \times {}^{J'_2} W)$. We call $[w_1, w_2, \cc, \cc']$ a $\car_{\cc'} \times \car_{\cc}$-stable piece of $G \times G$. 

Now we list some properties of the $\car_{\cc'} \times \car_{\cc}$-stable pieces. 

(1) The $\car_{\cc'} \times \car_{\cc}$-stable piece $[w_1, w_2, \cc, \cc']$ is a locally closed, smooth and irreducible subvariety of $G \times G$ of dimension equal to $\dim(G)+|I|+l(w_1)+l(w_2)+l(w_0^{J_1})+l(w_0^{J_2})$. See \cite[Theorem 2.2 (i)]{LY}. See also \cite[Theorem 2.6]{Sp3}. 

(2) $G \times G=\sqcup_{w_1 \in W^{J_1}, w_2 \in {}^{J'_2} W} [w_1, w_2, \cc, \cc']$. Lu and Yakimov \cite[2.2]{LY}  and Springer \cite[Theorem 2.6]{Sp3} gave two different proofs of this result. A different approach is sketched in \cite[Proposition 5.6]{H4}.

(3) Let $w_1 \in W^{J_1}$ and $w_2 \in {}^{J'_2} W$ and $\co=W_{c'} (w_1, w_2) W_c$. Then for any $(w'_1, w'_2) \in \co_{\min}$, $[\co, \cc, \cc']=\car_{\cc'} (B \dot w'_1 B, B \dot w'_2 B) \car_{\cc}$. See \cite[Proposition 5.3]{H4}. 

(4) Let $(w_1, w_2) \in ^{J'_1} W_1 \times W_2^{J_2}$. Define an automorphism $\th_{\s}: L_{I(w_1, w_2, c, c')} \rightarrow L_{I(w_1, w_2, c, c')}$ by $\th_{\s}(l)=\th_{\d} \i \bigl(w_2 \i \th_{\d} (w_1 l w_1 \i) w_2 \bigr)$. Then map $L_{I(w_1, w_2, c, c')} \rightarrow G_1 \times G_2$ defined by $l \rightarrow (w_1 l, w_2)$ induces a bijection between the $\th_{\s}$-twisted conjugacy classes on $L_{I(w_1,
w_2, c, c')}$ and the double cosets $\car_{\cc'} \backslash [w_1, w_2, \cc, \cc'] /\car_{\cc}$. See \cite[2.2]{LY} and \cite[Proposition 5.6 (2)]{H4}. 

(5) For any $(w_1, w_2) \in W \times W$, $\overline{\car_{\cc'} (B \dot w_1 B, B \dot w_2 B) \car_{\cc}}=\sqcup_{\co} [\co, \cc, \cc']$, where $\co$ runs over the distinguished double cosets in $W_{c'} \backslash (W \times W)/W_c$ that contains a minimal length element $(w'_1, w'_2)$ with $w'_1 \le w_1$ and $w'_2 \le w_2$. See \cite[Proposition 5.8]{H4}. A slightly more complicated description was obtained in \cite[Theorem 5.2]{LY}. 

In particular, 

(6) for any distinguished double coset $\co \in W_{c'} \backslash (W \times W)/W_c$, we have that $\overline{[\co, \cc, \cc']}=\sqcup_{\co' \le \co} [\co', \cc, \cc']$. See \cite[Corollary 5.9]{H4}.

\

Now we will come to the wonderful compactifications and the $P_K$-stable-piece decompositions on the compactifications.

From now on, unless otherwise stated, we assume that $G$ is adjoint and $\tilde G$ an algebraic group with identity component $G$. Let $G^1$ be a connected component of $\tilde G$. We fix an element $g_0 \in G^1$ with ${}^{g_0} B=B$ and ${}^{g_0} T=T$. If $G^1=G$, then we choose $g_0=1$ and $\d=id$. We denote by $\th_{\d}$ the conjugation of $g_0$ on $G$. Then $\th_{\d}$ gives automorphisms on $I$ and $W$. We denote these automorphisms by $\d$. 

\subsection{} We consider $G$ as a $G \times G$-variety by left and right translation. Let $\bar{G}$ be the wonderful compactification of $G$. This compactification was first constructed by De Concini and Procesi \cite{DP} when $\kk=\mathbf C$ and later generalized by Strickland \cite{Str} to arbitrary algebraically closed field $\kk$. It is known that $\bar{G}$ is an irreducible, smooth projective $(G \times G)$-variety with finitely many $G \times G$-orbits $Z_J$ indexed by the subsets $J$ of $I$. Here $Z_J$ is isomorphic to the quotient space $(G \times G) \times_{P^-_J \times P_J} G_J$ for the $P^-_J \times P_J$-action on $G \times G \times G_J$ defined by $(q, p) \cdot (g_1, g_2, z)=(g_1 q \i, g_2 p \i, \bar q z \bar p \i)$, where $\bar q$ is the image of $q$ under the projection $P^-_J \to G_J$ and $\bar p$ is the image of $p$ under the projection $P_J \to G_J$. Let $h_J$ be the image of $(1, 1, 1)$ in $Z_J$ under this isomorphism. 

\subsection{}\label{11} The wonderful compactification $\overline{G^1}$ of $G^1$ is the $(G \times G)$-variety which is isomorphic to $\bar{G}$ as a variety and where the $G \times G$-action is twisted by $(g, g') \mapsto (g, \th_{\d}(g'))$. The $G \times G$-orbits on $\overline{G^1}$ then coincides with the $G \times G$-orbits on $\bar{G}$. Let $Z_{J, \d}$ be the orbit coinciding with $Z_{\d(J)}$ and $h_{J, \d} \in Z_{J, \d}$ be the point identified with the base point $h_{\d(J)} \in Z_{\d(J)}$. Then $G^1$ is identified with the open $G \times G$-orbit $Z_{I, \d}$ via $g g_0 \mapsto (g, 1) \cdot h_{I, \d}$. Moreover, the isotropy subgroup of $h_{J, \d}$ in $G \times G$ is $$(U_{P^-_{\d(J)}} \times U_{P_J} Z(L_J)) (L_J)_{\d},$$ where $(L_J)_{\d}=\{(\th_{\d}(l), l); l \in L_J\}$. 

In other words, we have the following commuting diagram

\[\xymatrix{ G  \ar[r]^{\cdot g_0}  \ar@{^(->}[d] & G^1 \ar@{^(->}[d] \\ \bar{G} \ar[r]^{r} & \overline{G^1},}\] where $r \bigl( (g_1, g_2) \cdot h_{\d(J)} \bigr)=(g_1, \th_{\d} \i(g_2)) \cdot h_{J, \d}$.

For any subvariety $X \subset \overline{G^1}$, we denote by $\bar X$ its closure.

\subsection{}\label{p} For $J \subset I$, set $J_1=w_0 w_0^{\d(J)} \d(J)$ and $\d'=\d \i \circ \Ad(w_0 w_0^{\d(J)}) \i: W_{J_1} \to W_J$. Then $c=(J_1, J, \d')$ is an admissible triple on $W \times W$. Set $\th_{\d'}=\th_{\d} \i \circ \Ad(\dot w_0 \dot w_0^{\d(J)}) \i: L_{J_1} \to L_J$. Then $\cc=(J_1, J, \th_{\d'})$ is an admissible triple on $G \times G$. We may identify $(G \times G)/ \car_{\cc} (1, Z(L_J))$ with $Z_{J, \d}$ as $G \times G$-variety via $(g_1, g_2) \mapsto (g_1 \dot w_0 \dot w^{\d(J)}_0, g_2) \cdot h_{J, \d}$. 

Let $K \subset I$ and $\cc'=(K, K, id)$. Then each $\car_{\cc'} \times \car_{\cc}$-stable piece of $G \times G$ is stable under the right action of $\car_{\cc} (1, Z(L_J))$. For $w \in W^{\d(J)}$ and $v \in {}^K W$, set \[[J, w, v]_{K, \d}=[w, v, \cc, \cc']/\car_{\cc} (1, Z(L_J))=(P_K)_{\D} (B \dot w, B \dot v) \cdot h_{J, \d}.\]

We call $[J, w, v]_{K, \d}$ a $P_K$-stable piece on $\overline{G^1}$. In the case where $K=\emptyset$, a $P_K$-stable piece is just a $B \times B$-orbit and we simply write $[J, w, v]_{\d}$ for $[J, w, v]_{\emptyset, \d}$. In the case where $K=I$, a $P_K$-stable piece is just Lusztig's $G$-stable piece introduced in \cite[Section 12]{L3} and we simply write $Z_{J, w; \d}$ for $[J, w, 1]_{I, \d}$. 

The following properties follows easily from the properties of $\car_{\cc'} \times \car_{\cc}$-stable pieces that we listed in subsection \ref{cc}.

(1) $[J, w, v]_{K, \d}$ is an irreducible, locally closed subvariety of $\overline{G^1}$ of dimension $l(w_0)+|J|+l(v)-l(w)+l(w_0^K)$. 

(2) $\overline{G^1}=\sqcup_{J \subset I, w \in W^{\d(J)}, v \in {}^K W} [J, w, v]_{K, \d}$. 

(3) For any $J \subset I$, $x \in W^{\d(J)}$ and $y \in W$, $\overline{(P_K)_{\D} \cdot [J, x, y]_{\d}} \cap Z_{J, \d}=\sqcup [J, w, v]_{K, \d}$, where $(w, v)$ runs over all elements in $W^{\d(J)} \times {}^K W$ such that there exists $a \in W_K$ and $b \in W_J$ such that $a w \d(b) w_0^{\d(J)}  w_0 \le x w_0^{\d(J)} w_0$,$a v b \le y$ and $l(a w \d(b) w_0^{\d(J)}  w_0)+l(a v b)=l(w w_0^{\d(J)} w_0)+l(v)$.

(4) For $x \in W^{\d(J)}$ and $y \in W$ with $l(y)-l(x)=l(v)-l(w)$ and there exists $a \in W_K$ and $b \in W_J$ such that $x=a w \d(b)$ and $y=a v b$, then we have that $(P_K)_{\D} \cdot [J, x, y]_{\d}=[J, w, v]_{K, \d}$. 

The following explicit description of the closure of a $P_K$-stable piece in $\overline{G^1}$ was obtained in \cite[Theorem 7.6]{LY}, which generalized results on the $B \times B$-orbit closures in \cite[Proposition 2.4]{Sp1} and \cite[Proposition 6.3]{HT2} and the $G$-stable-piece closures in \cite[Theorem 4.5]{H3}. 

(5) $\overline{[J, w, v]_{K, \d}}$ is a union of $P_K$-stable pieces. Moreover, $[J', w', v']_{K, \d} \subset \overline{[J, w, v]_{K, \d}}$ if and only if $J' \subset J$ and there exists $x \in W_K$ and $y \in W_J$ such that $x w' \ge w \d(y)$ and $x v' \le v y$. 

From (1) and (5), we have the following useful consequence.

(6) For $J' \subset J$, $\dim(\overline{[J, w, v]_{K, \d}}) \cap Z_{J', \d}=\dim([J, w, v]_{K, \d})-|J|+|J'|$. 

We also need the following variation of subsection \ref{cc} (4),

(7) $[J, w, v]_{K, \d}=(P_K)_{\D} (L_{K_1} \dot w, \dot v) \cdot h_{J, \d}$, where $K_1=\max\{K' \subset K; w \i(K') \subset J, w \i(K')=\d(v \i(K'))\}$. 

Moreover, we have an explicit description of the semi-stable locus $\overline{G^1}^{ss}$ for the diagonal $G$-action on $\overline{G^1}$ in terms of $G$-stable pieces (see \cite{HS}). The case where $G^1=G$ was also studied by De Concini, Kannan and Maffei in \cite{DKM}. 

(8) $\overline{G^1}^{ss}=\sqcup_{J \subset I} Z_{J, 1; \d}$. 

\section{Closure of a parabolic subgroup in $\overline{G^1}^{ss}$}

For any $J \subset I$, set $J_{\d}=\max\{J_1 \subset J; \d(J_1)=J_1\}$. 

For any $K \subset I$ with $\d(K)=K$, we write $P_K^1=P_K g_0=N_{\tilde G} P_K \cap G^1$ and $G_K^1=L_K g_0/Z(L_K)$. Now we give an explicit description of $\overline{P_K^1} \cap \overline{G^1}^{ss}$ using $P_K$-stable pieces.

\begin{thm}\label{bp}
For $K \subset I$ with $\d(K)=K$,  we have that \[\overline{P_K} \cap \overline{G^1}^{\ss}=\sqcup_{J \subset I} \sqcup_{w \in {}^K W^J, w W_J \cap W^{\d} \neq \emptyset} [J, \d(w), w]_{K, \d}.\]
\end{thm}

Proof. By subsection \ref{p} (5), $$\overline{P^1_K} \cap Z_{J, 1; \d}=\sqcup_{w \in W^{\d(J)}, v \in {}^K W, x w \ge \d(x v) \text{ for some } x \in W_K} ([J, w, v]_{K, \d} \cap Z_{J, 1; \d}).$$

Let $w \in W^{\d(J)}, v \in {}^K W$ with $x w \ge \d(x v)$ for $x \in W_K$. Since $v \in {}^K W$, we have that \[l(w) \ge l(x w)-l(x) \ge l(x v)-l(x)=l(v).\] 

By subsection \ref{p} (3), $\overline{G_{\D} \cdot [J, w, v]_{K, \d}} \cap Z_{J, \d}=\overline{G_{\D} \cdot [J, w, v]_{\d}} \cap Z_{J, \d}$ is a union of $G$-stable pieces. If $[J, w, v]_{K, \d} \cap Z_{J, 1; \d} \neq \emptyset$, then by $Z_{J, 1; \d} \subset \overline{G_{\D} \cdot [J, w, v]_{\d}}$. Again by subsection \ref{p} (3), there exists $a \in W$ and $b \in W_J$ such that $a \d(b) w_0^{\d(J)} w_0 \le w w_0^{\d(J)} w_0$, $a b \le v$ and $l(a \d(b) w_0^{\d(J)} w_0)+l(a b)=l(w_0^{\d(J)} w_0)$. Therefore $a \d(b) w_0^{\d(J)} \ge w w_0^{\d(J)}$ and \begin{align*} l(a \d(b) w_0^{\d(J)}) & \ge l(w w_0^{\d(J)})=l(w)+l(w_0^{\d(J)}) \ge l(v)+l(w_0^{\d(J)}) \\ & \ge l(a b)+l(w_0^{\d(J)}).\end{align*}

Since $l(a \d(b) w_0^{\d(J)} w_0)+l(a b)=l(w_0^{\d(J)} w_0)$, we have that $$l(a \d(b) w_0^{\d(J)})=l(a b)+l(w_0^{\d(J)}).$$ Therefore, $x w=\d(x v)$, $a \d(b) w_0^{\d(J)}=w w_0^{\d(J)}$ and $a b=v$. So $w \d(b) \i=v b \i=a$ and $x v b\i=x w \d(b) \i=\d(x v) \d(b) \i=\d(x v b \i)$. 

We may write $x v b \i$ as $x v b \i=z_1 z_2$ for $z_1 \in W_K$ and $z_2 \in {}^K W$. Then $x v b \i=\d(x v b \i)=\d(z_1) \d(z_2)$ and $\d(z_1) \in W_K$, $\d(z_2) \in {}^K W$. Therefore $z_1=\d(z_1)$ and $z_2=\d(z_2)$. Write $z_2$ as $z_2=z_3 z_4$, where $z_3 \in W^J$ and $z_4 \in W_J$. Then $z_3 \in {}^K W^J$ and $x w \d(b) \i=x v b \i=z_1 z_3 z_4=z_1 \d(z_3) \d(z_4)$. By \cite[Corollary 2.5]{H4}, $(w, v)=(\d(z_3), z_3)$. 

Therefore $\overline{P^1_K} \cap Z_{J, 1; \d} \subset \sqcup_{J \subset I} \sqcup_{z \in {}^K W^J, z W_J \cap W^{\d} \neq \emptyset} [J, \d(z), z]_{K, \d}$. 

Now for $z \in {}^K W^J$ such that $z u=\d(z u)$ for some $u \in W_J$, we have that $G_{\D} \cdot [J, \d(z), z]_{K, \d}=G_{\D} \cdot [J, \d(z), z]_{\d}$. By subsection \ref{p} (4), $G_{\D} \cdot [J, \d(z), z]_{\d}=Z_{J, 1; \d}$. Hence $[J, \d(z), z]_{K, \d} \subset \overline{P^1_K} \cap Z_{J, 1; \d}$. The theorem is proved. \qed

\begin{lem}\label{7} Let $J, K \subset I$ with $\d(K)=K$. Then the map $w \mapsto \min(w W_J)$ gives a bijection $$\e: {}^K W^{J_{\d}} \cap W^{\d} \to \{x \in {}^K W^J, x W_J \cap W^{\d} \neq \emptyset\}.$$ Moreover, $\max\{K' \subset K; K'=\d(K'), \e(w) \i(K') \subset J\}=K \cap w(J_{\d})$.  
\end{lem}

Proof. If $w \in {}^K W^{J_{\d}} \cap W^{\d}$ and $x \in \min(w W_J)$. Then $x \in {}^K W^J$ and $w \in x W_J \cap W^{\d}$. So the map is well-defined. 

Now suppose that $x \in {}^K W^J$ with $x W_J \cap W^{\d} \neq \emptyset$. Let $y \in x W_J \cap W^{\d}$. Write $y$ as $y=a b$ for $a \in W_K$ and $b \in {}^K W$. Since $\d(K)=K$, we have that $\d(a) \in W_K$ and $\d(b) \in {}^K W$. Now $a b=y=\d(y)=\d(a) \d(b)$. So $b=\d(b)$. Since $b \in W_K x W_J \cap {}^K W$ and $x \in {}^K W^J$, we have that $b \in x W_J$. 

Write $b$ as $b=w c$ for $w \in {}^K W^{J_{\d}}$ and $c \in W_{J_{\d}}$. Then $w c=b=\d(b)=\d(w)\d(c)$ and $\d(w) \in {}^K W^{J_{\d}}$, $\d(c) \in W_{J_{\d}}$. Thus $w=\d(w) \in {}^K W^{J_{\d}} \cap W^{\d}$ and $\e(w)=x$. The map is surjective.

If $w_1, w_2 \in {}^K W^{J_{\d}}$ with $\e(w_1)=\e(w_2)$. Then $w_2=w_1 a$ for some $a \in W_J$. Thus $w_1 a=w_2=\d(w_2)=\d(w_1) \d(a)=w_1 \d(a)$ and $a=\d(a)$. Let $\supp(a)$ be the set of simple roots whose associated simple reflections appear in a reduced expression of $a$. Then $\supp(a)=\d(\supp(a)) \subset J$. Hence $\supp(a) \subset J_{\d}$ and $a \in W_{J_{\d}}$. Since $w_1, w_2 \in W^{J_{\d}}$, we have that $a=1$ and $w_1=w_2$. The map is injective. 

Let $w \in {}^K W^{J_{\d}} \cap W^{\d}$. Then $w=\e(w) a$ for some $a \in W_J \cap W^{J_{\d}}$. Let $K' \subset K$. If $a \i \e(w) \i(K')=w \i(K') \subset J_{\d} \subset J$, then $\e(w) \i(K') \subset \Phi_J$. Since $\e(w) \in W^J$, we must have that $\e(w) \i(K') \subset J$. Moreover, $\d(K \cap w J_{\d})=\d(K) \cap \d(w) \d(J_{\d})=K \cap w J_{\d}$. Hence $K \cap w (J_{\d}) \subset \max\{K' \subset K; \d(K')=K', \e(w) \i(K') \subset J\}$. On the other hand, assume that $K' \subset K$, $\d(K')=K$ and $\e(w) \i(K') \subset J$. Then for any $i \in K'$, $w \i(\a_i)=a \i \e(w) \i(\a_i)$ is a root in $\Phi_J$. Since $w \in {}^K W$, $w \i(\a_i)$ is a positive root in $\Phi_J$. Now $$\d \bigl( w\i \sum_{i \in K'} \a_i)=\d(w) \i \sum_{i \in \d(K')=K'} \a_i=w \i \sum_{i \in K'} \a_i.$$ Hence, $w \i(\a_i)$ is a positive root in $\Phi_{J_{\d}}$ for $i \in K'$. Notice that $w \in W^{J_{\d}}$. Thus $w \i(\a_i)$ is a simple root in $\Phi_{J_{\d}}$ for $i \in K'$ and $w \i(K') \subset J_{\d}$. \qed

\

Notice that for any $w \in {}^K W \cap W^{\d}$ and $J \subset I$, $\min(w W_{J_{\d}}) \in W^{\d}$ and \begin{align*} (P_K)_{\D} (B \dot w, B \dot w) \cdot h_{J, \d} &=(P_K)_{\D} \bigl(B \min(w W_{J_{\d}}), B \min(w W_{J_{\d}}) \bigr) \cdot h_{J, \d} \\ &=[J, \d(\min(w W_J)), \min(w W_J)]_{K, \d}. \end{align*}

By Proposition \ref{bp} and the previous lemma, we have other descriptions of $\overline{P^1_K} \cap \overline{G^1}^{ss}$ which are sometimes more convenient to use. 

\begin{thm}\label{pp} For $K \subset I$ with $\d(K)=K$, we have that \begin{align*} \overline{P_K^1} \cap \overline{G^1}^{ss} &=\sqcup_{J \subset I} \sqcup_{w \in {}^K W^{J^\d} \cap W^{\d}} (P_K)_{\D} (B \dot w, B \dot w) \cdot h_{J, \d} \\ &=\sqcup_{J \subset I} \cup_{w \in {}^K W \cap W^{\d}} (P_K)_{\D} (B \dot w, B \dot w) \cdot h_{J, \d} \\ &=\cup_{w \in {}^K W \cap W^{\d}} \sqcup_{J \subset I} (P_K)_{\D} (B \dot w, B \dot w) \cdot h_{J, \d}.\end{align*}
\end{thm}

\begin{thm}\label{smooth} For any $K \subset I$ with $\d(K)=K$, the variety $\overline{P_K^1} \cap \overline{G^1}^{ss}$ is smooth.
\end{thm}

The proof will be given in subsection \ref{pf}. The main idea of the proof is to find an open covering of $\overline{P_K^1} \cap \overline{G^1}^{ss}$ such that each open subvariety appeared in the covering is open in another smooth variety. 

\begin{lem} For any $K \subset I$ with $\d(K)=K$ and $w \in {}^K W$ with $\d(w)=w$, $\sqcup_{J \subset I} (B \dot w, B \dot w_0^K \dot w) \cdot h_{J, \d}$ is a locally closed subvariety of $\overline{G^1}$ isomorphic to an affine space of dimension $\dim(P_K)$. 

\end{lem}

Proof. Since $w \in {}^K W$ and $\d(w)=w$, we have that \begin{gather*} \th_{\d} \i ({}^{\dot w \i} U \cap U^-)={}^{\dot w \i} U \cap U^-={}^{\dot w \i} U_{P_K} \cap U^- \subset {}^{\dot w \i \dot w_0^K} U \cap U^-, \\ \th_{\d}({}^{\dot w \i \dot w_0^K} U \cap U)={}^{\dot w \i \dot w_0^K} U \cap U={}^{\dot w \i} U_{P_K} \cap U \subset {}^{\dot w \i} U \cap U. \end{gather*} 

For $J \subset I$, we have that \begin{align*} & ({}^{\dot w \i} B,  {}^{\dot w \i \dot w_0^K} B) \cdot h_{J, \d}= ({}^{\dot w \i} U, {}^{\dot w \i \dot w_0^K} B) \cdot h_{J, \d} \\ &= \bigl(({}^{\dot w \i} U \cap U) ({}^{\dot w \i} U \cap U^-), {}^{\dot w \i \dot w_0^K} B \bigr) \cdot h_{J, \d} \\ &=\bigl({}^{\dot w \i} U \cap U , {}^{\dot w \i \dot w_0^K} B \th_{\d} \i ({}^{\dot w \i} U \cap U^- \cap L_{\d(J)}) \bigr) \cdot h_{J, \d} \\ &=\bigl({}^{\dot w \i} U \cap U , {}^{\dot w \i \dot w_0^K} B \bigr) \cdot h_{J, \d} \\ &=\bigl({}^{\dot w \i} U \cap U , ({}^{\dot w \i \dot w_0^K} B \cap B^-) ({}^{\dot w \i \dot w_0^K} U \cap U) \bigr) \cdot h_{J, \d} \\ &=\bigl( ({}^{\dot w \i} U \cap U) \th_{\d} ({}^{\dot w \i \dot w_0^K} U \cap U \cap L_J), ({}^{\dot w \i \dot w_0^K} B \cap B^-) \bigr) \cdot h_{J, \d} \\ &=({}^{\dot w \i} U \cap U, {}^{\dot w \i \dot w_0^K} B \cap B^-) \cdot h_{J, \d}. \end{align*}

Set $X=\sqcup_{J \subset I} (1, T) \cdot h_{J, \d}$. Using the result of \cite[3.7 \& 3.8]{DS}, we see that $$(\dot w \i, \dot w \i  \dot w_0^K) \cdot \sqcup_{J \subset I} (B \dot w, B \dot w_0^K \dot w) h_{J, \d}=({}^{\dot w \i} U \cap U, {}^{\dot w \i \dot w_0^K} U \cap U^-) \cdot X$$ is a closed subvariety of $(U, U^-) \cdot X$ isomorphic to an affine space of dimension $\dim(P_K)$. 

Since $(U, U^-) \cdot X$ is open in $\overline{G^1}$, $(\dot w \i, \dot w \i \dot w_0^K) \cdot \sqcup_{J \subset I} (B \dot w, B \dot w_0^K \dot w) h_{J, \d}$ is locally closed in $\overline{G^1}$. \qed

\begin{lem}\label{1} For any $K \subset I$ with $\d(K)=K$ and $w \in {}^K W$, we have that $\sqcup_{J \subset I} (P_K \dot w, P_K \dot w) \cdot h_{J, \d}$ is smooth. 
\end{lem}

Proof. Set $X=\sqcup_{J \subset I} (B \dot w, B \dot w_0^K \dot w) \cdot h_{J, \d}$. Then $X$ is isomorphic to an affine space and \[\sqcup_{J \subset I} (P_K \dot w, P_K \dot w) \cdot h_{J, \d}=\cup_{p, q \in P_K} (p, q) \cdot X.\] So it suffices to prove that $X$ is open in $\sqcup_{J \subset I} (P_K \dot w, P_K \dot w) \cdot h_{J, \d}$. 

Suppose that $X$ is not open in $\sqcup_{J \subset I} (P_K \dot w, P_K \dot w) \cdot h_{J, \d}$. Notice that $X$ and \[\sqcup_{J \subset I} (P_K \dot w, P_K \dot w) \cdot h_{J, \d}=\sqcup_{J \subset I} \cup_{x, y \in W_K} (B \dot x \dot w, B \dot y \dot w) \cdot h_{J, \d}\] are unions of some $B \times B$-orbits. Thus there exists a $B \times B$-orbit $\co$ in $\sqcup_{J \subset I} (P_K \dot w, P_K \dot w) \cdot h_{J, \d}-X$ whose closure contains a $B \times B$-orbit $\co'$ in $X$. 

We may assume that $\co \subset Z_{J, \d}$ and $\co' \subset Z_{J', \d}$. Set $w'=\min(w W_J)$. Then $w' \in {}^K W^J$ and \begin{align*} \dim \bigl((B \dot w, B \dot w_0^K \dot w) \cdot h_{J, \d} \bigr) &=\dim([J, \d(w'), w_0^K w']_{\d})=l(w_0)+|J|+l(w_0^K) \\ &=\dim \bigl( (P_K \dot w, P_K \dot w) \cdot h_{J, \d} \bigr).\end{align*} Thus $(B \dot w, B \dot w_0^K \dot w) \cdot h_{J, \d}$ is open in $(P_K \dot w, P_K \dot w) \cdot h_{J, \d}$ and 
$\dim(\co)<\dim \bigl((B \dot w, B \dot w_0^K \dot w) \cdot h_{J, \d} \bigr)$. By subsection \ref{p} (6), \begin{align*} \dim(\overline{\co} \cap Z_{J'}) &<\dim \bigl((B \dot w, B \dot w_0^K \dot w) \cdot h_{J, \d} \bigr)-|J|+|J'| \\ &=\dim \bigl((B \dot w, B \dot w_0^K \dot w) \cdot h_{J', \d} \bigr)=\dim(\co').\end{align*} Therefore $\co'  \nsubseteq \overline{\co}$, which is a contradiction. \qed 

\begin{lem}\label{2} For any $K \subset I$ with $\d(J)=J$ and $w \in {}^K W$ with $\d(w)=w$, $\sqcup_{J \subset I} (P_K)_{\D} (B \dot w, B \dot w) \cdot h_{J, \d}$ is open in $\sqcup_{J \subset I} (P_K \dot w, P_K \dot w) \cdot h_{J, \d}$. 
\end{lem}

Proof. By definition, $(P_K)_{\D} (B \dot w, B \dot w) \cdot h_{J, \d}=[J, \d(\min(w W_J)), \min(w W_J)]_{K, \d}$. Thus $\sqcup_{J \subset I} (P_K)_{\D} (B \dot w, B \dot w) \cdot h_{J, \d}$ is a union of $P_K$-stable pieces.

Notice that $(P_K)_{\D} \subset P_K \times P_K$ and $B \times B \subset P_K \times P_K$. Thus for any $J \subset I$ and $x, v \in W$, either $$(P_K)_{\D} (B \dot x, B \dot v) \cdot h_{J, \d} \cap (P_K \dot w, P_K \dot w) \cdot h_{J, \d}=\emptyset$$ or $$(P_K)_{\D} (B \dot x, B \dot v) \cdot h_{J, \d} \subset (P_K \dot w, P_K \dot w) \cdot h_{J, \d}.$$ In other words, $\sqcup_{J \subset I} (P_K \dot w, P_K \dot w) \cdot h_{J, \d}$ is a union of $P_K$-stable pieces. 

Suppose that $\sqcup_{J \subset I} (P_K)_{\D} (B \dot w, B \dot w) \cdot h_{J, \d}$ is not open in $\sqcup_{J \subset I} (P_K \dot w, P_K \dot w) \cdot h_{J, \d}$. Then there exists a $P_K$-stable piece $\co$ in $\sqcup_{J \subset I} (P_K \dot w, P_K \dot w) \cdot h_{J, \d}-\sqcup_{J \subset I} (P_K)_{\D} (B \dot w, B \dot w) \cdot h_{J, \d}$ whose closure contains a $P_K$-stable piece $\co'$ in $\sqcup_{J \subset I} (P_K)_{\D} (B \dot w, B \dot w) \cdot h_{J, \d}$.

We may assume that $\co \in Z_{J, \d}$ and $\co' \in Z_{J', \d}$. Set $w'=\min(w W_J)$. By subsection \ref{p} (1), \begin{align*} \dim\bigl((P_K)_{\D} (B \dot w, B \dot w) \cdot h_{J, \d} \bigr) &=\dim([J, \d(w'), w']_{K, \d})=l(w_0)+|J|+l(w_0^K) \\ &=\dim\bigl((P_K \dot w, P_K \dot w) \cdot h_{J, \d} \bigr).\end{align*} Thus $[J, w, w]_{K, \d}$ is open in $(P_K \dot w, P_K \dot w) \cdot h_{J, \d}$ and $\dim(\co)< \dim([J, w', w']_{K, \d})$. By subsection \ref{p} (6), \begin{align*} \dim(\overline{\co} \cap Z_{J', \d}) &<\dim([J, w', w']_{K, \d})-|J|+|J'| \\ &=\dim \bigl((P_K)_{\D} (B \dot w, B \dot w) \cdot h_{J', \d} \bigr)=\dim(\co').\end{align*} Therefore $\co'  \nsubseteq \overline{\co}$, which is a contradiction. \qed 

\subsection{Proof of Theorem \ref{smooth}}\label{pf}

By Lemma \ref{1} and \ref{2}, for $w \in {}^K W \cap W^{\d}$, $\sqcup_{J \subset I} (P_K)_{\D} (B \dot w, B \dot w) \cdot h_{J, \d}$ is an open subvariety of a smooth variety. So $\sqcup_{J \subset I} (P_K)_{\D} (B \dot w, B \dot w) \cdot h_{J, \d}$ is smooth. By theorem \ref{pp}, $\sqcup_{J \subset I} (P_K)_{\D} (B \dot w, B \dot w) \cdot h_{J, \d} \subset \overline{P_K^1} \cap \overline{G^1}^{ss}$. Now it suffices to prove that $\sqcup_{J \subset I} (P_K)_{\D} (B \dot w, B \dot w) \cdot h_{J, \d}$ is open in $\overline{P^1_K} \cap \overline{G^1}^{ss}$ for any $w \in {}^K W \cap W^{\d}$. 

If this is not true, then there exists a $P_K$-stable piece $[J, \d(x), x]_{K, \d}$ in $\overline{P_K} \cap \overline{G^1}^{ss}-\sqcup_{J' \subset I} (P_K)_{\D} (B \dot w, B \dot w) \cdot h_{J', \d}$ whose closure contains a $P_K$-stable piece $(P_K)_{\D} (B \dot w, B \dot w) \cdot h_{J', \d}$. So $J' \subset J$. Set $w'=\min(w W_{J'})$. By subsection \ref{p} (5), there exists $a \in W_K$ and $b \in W_J$ such that $a \d(w') \ge \d(x) \d(b)$ and $a w' \le x b$. Then $$l(a)+l(w')=l(a \d(w')) \ge l(\d(x b))=l(x b) \ge l(a w')=l(a)+l(w').$$ Thus $a \d(w')=\d(x b)$ and $a w'=x b$. So $$\min(w' W_J)=\min(W_K a w' W_J)=\min(W_K x b W_J)=x$$ and $w \in w' W_{J'} \subset x W_J$. Assume that $w=x c$ for $c \in W_J$. Then $(B \dot w, B \dot w) \cdot h_{J, \d}=(B \dot{\d(x)} \dot{\d(c)}, B \dot x \dot c) \cdot h_{J, \d}=(B \dot{\d(x)}, B \dot x) \cdot h_{J, \d}$ and $[J, \d(x), x]_{K, \d}=(P_K)_{\D} (B \dot w, B \dot w) \cdot h_{J, \d}$, which is a contradiction. \qed

\section{A stratification on $\overline{G^1}^{ss}$}

First, we recall a stratification of $G^1$ introduced by Lusztig in \cite{L2}. 

\subsection{} An element $g \in G^1$ is called isolated if there is no proper parabolic subgroup $P$ of $G$ such that $h \in N_{\tilde G}(P)$ and $Z_G(h_s)^0 \subset P$, where $h_s$ is the semisimple part of $h$ (\cite[2.2]{L2}). Then the set of isolated elements is closed in $G^1$ (\cite[Lemma 2.8]{L2}) and the action of $Z(G) \times G$ on $G^1$ defined by $(z, g) \cdot g'=g z g' g \i$ leaves stable the set of isolated elements in $G^1$ and there are finitely many orbits there (\cite[Lemma 2.7]{L2}). These orbits are called isolated strata of $G^1$ (\cite[3.3]{L2}). 

\subsection{}\label{sm1} Let $P$ be a parabolic subgroup of $G$, $L$ be a Levi subgroup of $P$ and $S$ be an isolated stratum of $N_{\tilde G} (L) \cap G^1$ such that $S \subset N_{\tilde G}(P)$. Set $S^*=\{g \in S; Z_G(g_s)^0 \subset L\}$ and $Y_{L, S}=\sqcup_{g \in G} g S^* g \i$. We call $Y_{L, S}$ a stratum of $G^1$. It is known that $Y_{L, S}$ is smooth (\cite[3.17]{L2}) and $Y_{L, S}$ (for various $(L, S)$) form a stratification of $G^1$ (\cite[Proposition 3.12 \& 3.15]{L2}). 

Moreover, let $S'$ be the closure of $S$ in $N_{\tilde G} (L) \cap G^1$ and $G \times_P (S' U_P)$ be the quotient space of $G \times (S' U_P)$ under the $P$-action defined by $p (g, z)=(g p \i, p z p \i)$. Then the proper map $f: G \times_P (S' U_P) \to \overline{Y_{L, S}}$ defined by $(g, z) \mapsto g z g \i$ is a small map. See the proof of \cite[Proposition 5.7]{L2}. 

\

Now we generalize the definition of strata to $\overline{G^1}^{ss}$. 

\subsection{}\label{4} By \cite[Proposition 1.10]{H3}, the map $(g, z) \mapsto (g, g) \cdot z$ gives an isomorphism $G \times_{P_{J_{\d}}} (P_{J_{\d}}, P_{J_{\d}}) \cdot h_{J, \d} \to Z_{J, 1; \d}$. 

Notice that the map $(g, z) \mapsto (g, 1) z$ gives an isomorphism from $U_{P_{J_{\d}}} \times (L_{J_{\d}}, 1) \cdot h_{J, \d}$ to $(P_{J_{\d}}, P_{J_{\d}}) \cdot h_{J, \d}$ and the action of $U_{P_{J_{\d}}}$ on $(L_{J_{\d}}, 1) \cdot h_{J, \d}$ defined by $(g, z) \mapsto (1, g) \cdot z$ is trivial. Then $(g, z) \mapsto (g, g) \cdot z$ gives an isomorphism \[\tag{a} U_{P_{J_{\d}}} \times (L_{J_{\d}}, 1) \cdot h_{J, \d} \cong (P_{J_{\d}}, P_{J_{\d}}) \cdot h_{J, \d}.\]

Therefore \begin{align*}  P_{J_{\d}} \times_{L_{J_{\d}}} (L_{J_{\d}}, 1) \cdot h_{J, \d} & \cong  (U_{P_{J_{\d}}} \times L_{J_{\d}}) \times_{L_{J_{\d}}} (L_{J_{\d}}, 1) \cdot h_{J, \d} \\ & \cong U_{P_{J_{\d}}} \times (L_{J_{\d}} \times_{L_{J_{\d}}} (L_{J_{\d}}, 1) \cdot h_{J, \d}) \\ & \cong (P_{J_{\d}}, P_{J_{\d}}) \cdot h_{J, \d}, \end{align*} where $P_{J_{\d}} \times_{L_{J_{\d}}} (L_{J_{\d}}, 1) \cdot h_{J, \d}$ is the quotient space for the $L_{J_{\d}}$ action on $P_{J_{\d}} \times (L_{J_{\d}}, 1) \cdot h_{J, \d}$ defined by $l \cdot (p, z)=(p l \i, (l, l) \cdot z)$. 

Thus $Z_{J, 1; \d}$ is isomorphic to $G \times_{P_{J_{\d}}} (P_{J_{\d}} \times_{L_{J_{\d}}} (L_{J_{\d}}, 1) \cdot h_{J, \d}) \cong G \times_{L_{J_{\d}}} (L_{J_{\d}}, 1) \cdot h_{J, \d}$ via $(g, z) \mapsto (g, g) \cdot z$. 

We may also identify $(L_{J_{\d}}, 1) \cdot h_{J, \d}$ with $L_{J_{\d}} g_0/Z(L_J)$. Therefore we have an isomorphism $$i_J: G \times_{L_{J_{\d}}} L_{J_{\d}} g_0/Z(L_J) \cong Z_{J, 1; \d}$$ via $(g, l g_0) \mapsto (g l, g) \cdot h_{J, \d}$. 

Notice that we have a stratification $G \times_{L_{J_{\d}}} L_{J_{\d}} g_0=\sqcup\, G \times_{L_{J_{\d}}} Y$, where $Y$ runs over strata of $L_{J_{\d}} g_0$. Moreover, each stratum $Y$ of $L_{J_{\d}} g_0$ is stable under the action of $Z(L_{J_{\d}}) \supset Z(L_J)$. Then \[\tag{b} \overline{G^1}^{ss}=\sqcup_{J \subset I} \sqcup_{Y \text{ is a  stratum of } L_{J_{\d}} g_0} i_J(G \times_{L_{J_{\d}}} Y/Z(L_J))\] is a decomposition of $\overline{G^1}^{ss}$.
We will see later that $(b)$ is in fact a stratification. For any $J \subset I$ and stratum $Y$ of $L_{J_{\d}} g_0$, we call $i_J(G \times_{L_{J_{\d}}} Y/Z(L_J))$ a stratum of $\overline{G^1}^{ss}$. 

We may define a decomposition for $\overline{G^1}$ in the same way. But it is very hard to give an explicit description of the closure of any subvariety appeared in the decomposition. However, \cite[Theorem 4.3]{H1}, \cite[Theorem 7.4]{HT1} and \cite[Theorem 4.5]{H3} give some evidence that this decomposition for $\overline{G^1}$ may still be a stratification. 

\subsection{}\label{abs} In this subsection, we assume that $G^1=G$. It is known \cite{DP} that the map $(g, g', z) \mapsto (g, g') \cdot z$ gives an isomorphism $$(G \times G) \times_{P_J^- \times P_J} \overline{G_J} \cong \overline{Z_J}.$$ Notice that any  element in $\overline{Z_J} \cap \overline{G}^{ss}$ is of the form $(g l, g) \cdot h_K$ for some $K \subset J$, $g \in G$ and $l \in L_K$ and any element in $\overline{G_J}^{ss}$ is of the form $(g' l', g') \cdot h_K$ for some $K \subset J$, $g' \in L_J$ and $l' \in L_K$. Therefore $\overline{Z_J} \cap \overline{G}^{ss}=G_{\D} \cdot \overline{G_J}^{ss}$. 

The morphism $(G \times G) \times_{P_J^- \times P_J} \overline{G_J} \to G/P_J^- \times G/P_J$, $(g, g', z) \mapsto (g P_J^-, g' P_J)$ sends $\overline{Z_J} \cap \overline{G}^{ss}$ to the open $G_{\D}$ orbit $\co$ in $G/P_J^- \times G/P_J$. It is easy to see that $\co \cong G/L_J$. Since each fiber of the $G$-equivariant morphism $\overline{Z_J} \cap \overline{G}^{ss} \to G/L_J$ is isomorphic to $\overline{G_J}^{ss}$, by \cite[Page 26, Lemma 4]{Sl}, we have that $$\overline{Z_J} \cap \overline{G}^{ss} \cong G \times_{L_J} \overline{G_J}^{ss}.$$ Here $G \times_{L_J} \overline{G_J}^{ss}$ is the quotient space for the $L_J$-action on $G \times \overline{G_J}^{ss}$ defined by $l \cdot (g, z)=(g l \i, (l, l) \cdot z)$. This isomorphism extends the isomorphism $Z_{J, 1; id}=Z_J \cap \overline{G}^{ss} \cong G \times_{L_J} G_J$ in the previous subsection. 

\begin{lem} Let $T_0=\{t \th_{\d}(t \i); t \in T\}$. Let $J, K \subset I$ with $\d(K)=K$ and $w \in W^{\d}$. Then \[\overline{T_0 Z(L_K)} \cap (\dot w T, \dot w) \cdot h_J=\begin{cases} (T_0 Z(L_K) \dot w, \dot w) \cdot h_J, & \text{ if } w \i \Phi_K \subset \Phi_J; \\ \emptyset, & \text{ otherwise}. \end{cases}\]
\end{lem}

Proof. Let $X=\sqcup_{D \subset I} (T, 1) \cdot h_D$. Then for any positive root $\a$, the morphism $T \to k$ defined by $t \mapsto \a(t)$ extends in a unique way to a morphism from $X \to k$, which we denote by $\tilde \a$. It is easy to see that \[\tag{a} \tilde \a_i((t, 1) \cdot h_J)=\begin{cases} \a_i(t), & \text{ if } i \in J, \\ 0, & \text{ if } i \notin J. \end{cases}\] 

By definition, $T_0 Z(L_K)=\{t \in T; \prod_{i \in \co} \a_i(t)=1, \forall \text{ $\d$-orbit } \co \text{ of } K\}$. So $(\dot w \i, \dot w \i) \cdot T_0 Z(L_K)=\{t \in T; \prod_{i \in \co} w \i \a_i(t)=1, \forall \text{ $\d$-orbit } \co \text{ of } K\}$. For any root $\a$, set \[\sgn(\a)=\begin{cases} 1, & \text{ if } \a>0; \\ -1, & \text{ if } \a<0. \end{cases} \] Notice that $\d(w \i \a_i)=\d(w) \i \a_{\d(i)}=w \i \a_{\d(i)}$. Thus for any $\d$-orbit $\co$ of $K$, either $w \i \a_i>0$ for all $i \in \co$ or $w \i \a_i<0$ for all $i \in \co$. So we may write $\sgn(w \i \co)$ for $\sgn(w \i \a_i)$, where $i \in \co$. Now $\bigl(\prod_{i \in \co} \widetilde{w \i(\a_i)} z \bigr)^{\sgn(w \i \co)}$ is a well-defined morphism from $X$ to $k$ and $$\overline{(\dot w \i, \dot w \i) \cdot T_0 Z(L_K)}=\{z \in X; \prod_{i \in \co} \widetilde{w \i(\a_i)}(z)^{\sgn(w \i \co)}=1\}.$$

By (a), if $w \i (\Phi_K) \nsubseteq \Phi_J$, then $\prod_{i \in K} \widetilde{w \i(\a_i)}(z)^{\sgn(w\i(\a_i)}=0$ for all $z \in (T, 1) \cdot h_J$ and $\overline{(\dot w \i, \dot w \i) \cdot T_0 Z(L_K)} \cap (T, 1) \cdot h_J=\emptyset$. On the other hand, if $w \i(\Phi_K) \subset \Phi_J$, then for any $z=(t, 1) \cdot h_J$ and $i \in K$, $\widetilde{w \i(\a_i)}(z)=w \i(\a_i)(t)$. Therefore $\overline{(\dot w \i, \dot w \i) \cdot T_0 Z(L_K)} \cap (T, 1) \cdot h_J=\{(t, 1) \cdot h_J; t \in (\dot w \i, \dot w \i) \cdot T_0 Z(L_K)\}$. The lemma is proved. \qed

\

Notice that $\th_{\d} (T_0)=T_0$ and $\th_{\d} Z(L_K)=Z(L_K)$ for $K \subset I$ with $\d(K)=K$. By the identification of $\bar{G}$ with $\overline{G^1}$ in subsection \ref{11}, we have the following variation of the previous lemma.

\begin{lem} Let $J, K \subset I$ with $\d(K)=K$ and $w \in W^{\d}$. Then $$\overline{T_0 Z(L_K) g_0} \cap (\dot w T, \dot w) \cdot h_{J, \d}=\begin{cases} (T_0 Z(L_K) \th_{\d} (\dot w), \dot w) \cdot h_{J, \d}, & \text{ if } w \i \Phi_K \subset \Phi_J; \\ \emptyset, & \text{ otherwise}. \end{cases}$$
\end{lem}

\begin{thm}\label{xxx} Let $K \subset I$ with $\d(K)=K$ and $S$ be an isolated stratum of $L_K g_0$. Let $S'$ be the closure of $S$ in $L_K g_0$. Then \begin{align*} \overline{U_{P_K} S} \cap \overline{G^1}^{ss} &=\sqcup_{J \subset I, w \in {}^K W^{J_{\d}} \cap W^{\d}, w \i(K) \subset J_{\d}} (U_{P_K} S' \dot w g_0 \i, U_{P_K} \dot w) \cdot h_{J, \d}. \end{align*} 
\end{thm}

Proof. Since $U_{P_K} S' \subset P_K^1$, then $$\overline{U_{P_K} S'} \cap \overline{G^1}^{ss} \subset \overline{P_K^1} \cap \overline{G^1}^{ss}=\sqcup_{J \subset I} \sqcup_{w \in {}^K W^{J_{\d}} \cap W^{\d}} (P_K)_{\D} (B \dot w, B \dot w) \cdot h_{J, \d}.$$ 

Since $S$ is stable under the conjugation action of $L_K$, there exists $s \in (B \cap L_K) g_0$ such that $s \in S$. We may write $s$ as $s=u t g_0$ for some $t \in T$ and $u \in U \cap L_K$. Then $$S=\{l u t Z(L_K) g_0 l \i; l \in L_K\} \subset (L_K)_{\D} \bigl((U \cap L_J) t T_0 Z(L_K) g_0 \bigr).$$ 

It is known that $(U_{P_K})_{\D} \cdot S$ is dense in $U_{P_K} S'$. Now consider the proper map $P_K \times_B \overline{U t T_0 Z(L_K) g_0} \to \overline{G^1}$ defined by $(g, z) \mapsto (g, g) \cdot z$. Since $(P_K)_{\D} (U t T_0 Z(L_K) g_0)=(U_{P_K})_{\D} (L_K)_{\D} \cdot (U t T_0 Z(L_K) g_0) \supset (U_{P_K})_{\D} \cdot S$, then \[\tag{a} \overline{U_{P_K} S'} \subset (P_K)_{\D} \overline{U t T_0 Z(L_K) g_0}.\] 

Let $J \subset I$ and $w \in {}^K W^{J_{\d}} \cap W^{\d}$. By definition, \begin{align*} (P_K)_{\D} (B \dot w, B \dot w) \cdot h_{J, \d} & \subset \cup_{x \in W_K} (B \dot x B \dot w, B \dot x B \dot w) \cdot h_{J, \d} \\ &=\cup_{x \in W_K} (B \dot x B \dot w, B \dot x \dot w) \cdot h_{J, \d}  \\ & \subset \cup_{x \in W_K} \bigl( (B \dot x \dot w, B \dot x \dot w) \cdot h_{J, \d} \cup \cup_{y<x w} (B \dot y, B \dot x \dot w) \cdot h_{J, \d} \bigr).\end{align*} 

If $(P_K)_{\D} (B \dot w, B \dot w) \cdot h_{J, \d} \cap \overline{U_{P_K} S'} \neq \emptyset$, by (a) we have that \begin{align*} & \cup_{x \in W_K^{\d}} (B \dot x \dot w, B \dot x \dot w) \cdot h_{J, \d} =(P_K)_{\D} (B \dot w, B \dot w) \cdot h_{J, \d} \cap \overline{B^1}  \\ & \supset (P_K)_{\D} (B \dot w, B \dot w) \cdot h_{J, \d} \cap \overline{U t T_0 Z(L_K) g_0} \neq \emptyset. \end{align*} 

Therefore $\overline{U t T_0 Z(L_K) g_0} \cap (B \dot x \dot w, B \dot x \dot w) \cdot h_{J, \d} \neq \emptyset$ for some $x \in W_K^{\d}$. Set $w'=x w y$ and $X=\cup_{J' \subset I} (B \dot x \dot w, B \dot x \dot w) \cdot h_{J', \d}$. Then the map $(u, u', z) \mapsto (u \dot x \dot w, u' \dot x \dot w) \cdot z$ defines an isomorphism $$(U \cap {}^{\dot x \dot w} U) \times (U \cap {}^{\dot x \dot w} U^-) \times \cup_{J' \subset I} (T, 1) \cdot h_{J', \d} \to X.$$

Notice that $U t T_0 Z(L_K) g_0 \subset X$. Then $\overline{U t T_0 Z(L_K) g_0} \cap X$ is the closure of $U t T_0 Z(L_K) g_0$ in $X$. Hence $$\overline{U t T_0 Z(L_K) g_0} \cap X=\bigl((U \cap {}^{\dot x \dot w} U) \dot x \dot w, (U \cap {}^{\dot x \dot w} U^-) \dot x \dot w \bigr) \cdot X',$$ where $X'=\overline{^{(\dot x \dot w) \i} T_0 Z(L_K) g_0} \cap \cup_{J' \subset I} (T, 1) \cdot h_{J', \d}$. Since $\overline{U t T_0 Z(L_K) g_0} \cap (B \dot x \dot w, B \dot x \dot w) \cdot h_{J, \d} \neq \emptyset$, then $X' \cap (T, 1) \cdot h_{J, \d} \neq \emptyset$. By the previous lemma, $w \i \Phi_K=w \i x \i \Phi_K \subset \Phi_J$. Since $w=\d(w)$ and $K=\d(K)$, we have that $\d(w \i \Phi_K)=w \i \Phi_K \subset \Phi_{J_{\d}}$. Notice that $w \in {}^K W^{J_{\d}}$. Then $w \i(K) \subset J_{\d}$. 

On the other hand, suppose that $w \in {}^K W \cap W^{\d}$ with $w \i(K) \subset I$. Set $$Y=\sqcup_{w \i(K) \subset D \subset I} (P_K)_{\D} (B \dot w, B \dot w) \cdot h_{D, \d}.$$ If $w \i(K) \subset D$, then $g_0 L_{w \i(K)} g_0 \i=L_{w \i(K)}$ and by \cite[Lemma 7.3]{St}, $(L_{w \i(K)})_{\D} \cdot \bigl((B \cap L_{w \i(K)}) g_0 \bigr)=L_{w \i(K)} g_0$. Hence \begin{align*} & (P_K)_{\D} (B \dot w, B \dot w) \cdot h_{D, \d}=(L_K)_{\D} (U_{P_K} \dot w, U_{P_K} \dot w (B \cap L_{w \i(K)}) \cdot h_{D, \d} \\ &=(U_{P_K} \dot w, U_{P_K} \dot w) (L_{w \i(K)})_{\D} (1, B \cap L_{w \i(K)}) \cdot h_{D, \d} \\ &=(U_{P_K} \dot w, U_{P_K} \dot w) (1, L_{w \i(K)}) \cdot h_{D, \d} \\ &=(U_{P_K} \dot w, U_{P_K} \dot w) (L_{w \i(K)}, L_{w \i(K)}) \cdot h_{D, \d}=(P_K \dot w, P_K \dot w) \cdot h_{D, \d}.\end{align*} Therefore $Y=\sqcup_{w \i(K) \subset D \subset I} (P_K \dot w, P_K \dot w) \cdot h_{D, \d}$. Since $w \i(K) \subset I$, $U_{P_K} \cap {}^{\dot w} U=U_{P_K} \cap {}^{\dot w} U_{P_{w \i(K)}}$ and $U_{P_K} \cap {}^{\dot w} U^-=U_{P_K} \cap {}^{\dot w} U_{P^-_{w \i(K)}}$. It is easy to see that the map $(u, u', z) \mapsto (u, u') \cdot z$ defines an isomorphism \[\tag{b} (U_{P_K} \cap {}^{\dot w} U) \times (U_{P_K} \cap {}^{\dot w} U^-) \times \sqcup_{w \i(K) \subset D \subset I} (L_K \dot w, L_K \dot w) \cdot h_{D, \d} \to Y.\]

By the similar argument as we did above, one can show that the closure of $S'=(S' g_0 \i \th_{\d}(\dot w), \dot w) \cdot h_{I, \d}$ in $\sqcup_{w \i(K) \subset D \subset I} (L_K \dot w, L_K \dot w) \cdot h_{D, \d}$ is $\sqcup_{w \i(K) \subset D \subset I} (S' g_0 \i \th_{\d}(\dot w), \dot w) \cdot h_{D, \d}$. 

Hence the closure of $U_{P_K} S'=(U_{P_K} \cap {}^{\dot w} U, U_{P_K} \cap {}^{\dot w} U^-) \cdot (S' g_0 \i \th_{\d}(\dot w), \dot w) \cdot h_{I, \d}$ in $Y$ is \begin{align*} & (U_{P_K} \cap {}^{\dot w} U, U_{P_K} \cap {}^{\dot w} U^-) \cdot \sqcup_{w \i(K) \subset D \subset I} (S' g_0 \i \th_{\d}(\dot w), \dot w) \cdot h_{D, \d} \\ &=\sqcup_{w \i(K) \subset D \subset I} (U_{P_K} S' g_0 \i \th_{\d}(\dot w), U_{P_K} \dot w) \cdot h_{D, \d} \\ &=\sqcup_{w \i(K) \subset D \subset I} (U_{P_K} S' \dot w g_0 \i, U_{P_K} \dot w) \cdot h_{D, \d}. \end{align*}

The theorem is proved. \qed

\begin{cor}
Let $K \subset I$ and $S$ be an isolated stratum of $L_K$. Let $S'$ be the closure of $S$ in $L_K$. Then for any $J \subset I$, $$\overline{U_{P_K} S} \cap \overline{G}^{ss} \cap \overline{Z_J}=\sqcup_{w \in {}^K W^J, w \i(K) \subset J} (U_{P_K} \dot w, U_{P_K} \dot w) \cdot X_w,$$ where $X_w$ is the closure of $(U_{P_{w \i(K)}} \cap L_J) {}^{\dot w \i} S/Z(L_J)$ in $\overline{G_J}^{ss}$. 
\end{cor}

\subsection{}\label{11x} Let $J, K \subset I$ with $\d(K)=K$ and $w \in {}^K W^{J_{\d}} \cap W^{\d}$. Let $S \subset L_{J_{\d}}$ be a subvariety. Since $Z_{J, 1; \d} \cong G \times_{P_{J_{\d}}} (P_{J_{\d}}, P_{J_{\d}}) \cdot h_{J, \d}$, we have a projection map $Z_{J, 1; \d} \to G/P_{J_{\d}}$. Restricting the projection map to $({}^{\dot w \i} P_K)_{\D} (S, 1) \cdot h_{J, \d} \subset Z_{J, 1; \d}$, we obtain a morphism $$({}^{\dot w \i} P_K)_{\D} (S, 1) \cdot h_{J, \d} \to {}^{\dot w \i} P_K/{}^{\dot w \i} P_K \cap P_{J_{\d}}.$$ By \cite[page 26, lemma 4]{Sl}, \[({}^{\dot w \i} P_K)_{\D} (S, 1) \cdot h_{J, \d} \cong {}^{\dot w \i} P_K \times_{{}^{\dot w \i} P_K \cap P_{J_{\d}}} ({}^{\dot w \i} P_K \cap P_{J_{\d}})_{\D} (S, 1) \cdot h_{J, \d}.\]

Notice that ${}^{\dot w \i} P_K \cap P_{J_{\d}} \cong {}^{\dot w \i} P_K \cap U_{P_{J_{\d}}} \times {}^{\dot w \i} P_K \cap L_{J_{\d}}$ and $${}^{\dot w \i} P_K \cap L_{J_{\d}}=L_{K'} ({}^{\dot w \i} U_{P_K} \cap L_{J_{\d}})=L_{K'} (B \cap L_{J_{\d}})=P_{K'} \cap L_{J_{\d}},$$ where $K'=w \i K \cap J_{\d}$. Therefore, $$({}^{\dot w \i}  P_K \cap P_{J_{\d}})_{\D} (S, 1) \cdot h_{J, \d}=({}^{\dot w \i} P_K \cap U_{P_{J_{\d}}})_{\D} (P_{K'} \cap L_{J_{\d}})_{\D} (S, 1) \cdot h_{J, \d}.$$

Set $X=(P_{K'} \cap L_{J_{\d}})_{\D} (S, 1) \cdot h_{J, \d}$. By subsection \ref{4} (a), the map $(g, z) \mapsto (g, g) \cdot z$ gives an isomorphism \[({}^{\dot w \i} P_K \cap P_{J_{\d}})  \times_{P_{K'} \cap L_{J_{\d}}} X  \cong ({}^{\dot w \i} P_K \cap U_{P_{J_{\d}}})  \times X \cong ({}^{\dot w \i} P_K \cap U_{P_{J_{\d}}})_{\D} \cdot X. \]

Therefore, we have that \begin{align*}\tag{a} & G \times_{{}^{\dot w \i} P_K} ({}^{\dot w \i} P_K)_{\D} \cdot X \\ & G \times_{{}^{\dot w \i} P_K} \bigl({}^{\dot w \i} P_K \times_{{}^{\dot w \i} P_K \cap P_{J_{\d}}} ({}^{\dot w \i} P_K \cap P_{J_{\d}})_{\D} \cdot X\bigr) \\ & \cong G \times_{{}^{\dot w \i} P_K \cap P_{J_{\d}}} \bigl({}^{\dot w \i} P_K \cap P_{J_{\d}})_{\D} \cdot X \bigr) \\ & \cong G \times_{{}^{\dot w \i} P_K \cap P_{J_{\d}}} \bigl(({}^{\dot w \i} P_K \cap P_{J_{\d}})  \times_{P_{K'} \cap L_{J_{\d}}} X \bigr) \\ & \cong G \times_{P_{K'} \cap L_{J_{\d}}} X \cong G \times_{L_{J_{\d}}} (L_{J_{\d}} \times_{P_{K'} \cap L_{J_{\d}}} X). \end{align*}

Similarly, we may identify $G_{\D} \cdot X$ with $G \times_{L_{J_{\d}}} (L_{J_{\d}})_{\D} \cdot X$ and under these identifications, the map $(g, z) \mapsto (g, g) \cdot z$ from $G \times_{{}^{\dot w \i} P_K} ({}^{\dot w \i} P_K)_{\D} \cdot X$ to $G_{\D} \cdot X$ is induced from the map $$G \times (L_{J_{\d}} \times_{P_{K'} \cap L_{J_{\d}}} X) \to G \times ((L_{J_{\d}})_{\D} \cdot X),$$ defined by $(g, l, z) \mapsto (g, (l, l) \cdot z)$. 

\subsection{}\label{bbbbb} Let $J, K \subset I$ with $\d(K)=K$ and $w \in {}^K W^{J_{\d}} \cap W^{\d}$ with $w \i(K) \subset J_{\d}$. Let $S' \subset L_K g_0$ be the closure of an isolated stratum.  We have that ${}^{\dot w \i} P_K={}^{\dot w \i} U_{P_K}  {}^{\dot w \i} L_K=({}^{\dot w \i} U_{P_K} \cap U_{P_{J_{\d}}}) ({}^{\dot w \i} P_K \cap L_{J_{\d}}) ({}^{\dot w \i} U_{P_K} \cap U_{P^-_{J_{\d}}})$. Since $w \in W^{J_{\d}}$ and $w \i(K) \subset J_{\d}$, we have that ${}^{\dot w \i} U_{P_K} \cap L_{J_{\d}}=U_{P_{w \i(K)}} \cap L_{J_{\d}}$. Then \begin{align*} & (\dot w \i, \dot w \i) (P_K)_{\D} (U_{P_K} S' \dot w g_0 \i, U_{P_K} \dot w) \cdot h_{J, \d} \\ &=({}^{\dot w \i} P_K)_{\D} \bigl(({}^{\dot w \i} U_{P_K} \cap L_{J_{\d}}) ({}^{\dot w \i} U_{P_K} \cap U_{P^-_{J_{\d}}}) {}^{\dot w \i} S' g_0 \i, ({}^{\dot w \i} U_{P_K} \cap U_{P_{J_{\d}}}) \bigr) \cdot h_{J, \d} \\ &=({}^{\dot w \i} P_K)_{\D} (({}^{\dot w \i} U_{P_K} \cap L_{J_{\d}}) {}^{\dot w \i} S' g_0 \i, 1) \cdot h_{J, \d} \\ &=({}^{\dot w \i} P_K)_{\D} ((U_{P_{w \i K}} \cap L_{J_{\d}}) {}^{\dot w \i} S' g_0 \i, 1) \cdot h_{J, \d}.\end{align*}

The map $f: G \times_{P_K} (P_K)_{\D} (U_{P_K} S' \dot w g_0 \i, U_{P_K} \dot w) \cdot h_{J, \d} \to G \times_{{}^{\dot w \i} P_K}  ({}^{\dot w \i} P_K)_{\D} (U_{P_{w \i K}} \cap L_{J_{\d}}) {}^{\dot w \i} S' g_0 \i, 1) \cdot h_{J, \d}$ defined by $(g, z) \mapsto (g \dot w, (\dot w \i, \dot w \i) z)$ is an isomorphism. Moreover, \[\tag{*} \pi=\pi' \circ f,\] where \begin{gather*} \pi: G \times_{P_K} (P_K)_{\D} (U_{P_K} S' \dot w g_0 \i, U_{P_K} \dot w) \cdot h_{J, \d} \to Z_{J, 1; \d}, \\ \pi':  G \times_{{}^{\dot w \i} P_K}  ({}^{\dot w \i} P_K)_{\D} (U_{P_{w \i K}} \cap L_{J_{\d}}) {}^{\dot w \i} S' g_0 \i, 1) \cdot h_{J, \d} \to Z_{J, 1; \d}, \end{gather*} are induced from the map $G \times Z_{J, 1; \d} \to Z_{J, 1; \d}$ defined by $(g, z) \mapsto (g, g) \cdot z$.

As in the previous subsection, the map $\pi'$ is induced from the map $G \times (L_{J_{\d}} \times_{P_{w \i(K)} \cap L_{J_{\d}}} X) \to G \times ((L_{J_{\d}})_{\D} \cdot X)$ defined by $(g, l, z) \mapsto (g, (l, l) \cdot z)$, here $$X=(U_{P_{w \i K}} \cap L_{J_{\d}}) {}^{\dot w \i} S' g_0 \i, 1) \cdot h_{J, \d} \cong (U_{P_{w \i K}} \cap L_{J_{\d}}) {}^{\dot w \i} S'/Z(L_{J_{\d}}).$$ Since ${}^{\dot w \i} S'/Z(L_{J_{\d}})$ is the closure of an isolated stratum in $L_{P_{w \i K}} g_0/Z(L_{J_{\d}})$, by subsection \ref{sm1}, $G_{\D} \cdot X$ is a union of strata in $Z_{J, 1; \d}$ and the map $\pi'$ is a small map. 

As a summary, we have the following result.

\begin{thm}\label{sm2} Let $K \subset I$ with $\d(K)=K$ and $S$ be an isolated stratum of $L_K g_0$. Then the proper map $$G \times_{P_K} (\overline{U_{P_K} S} \cap \overline{G^1}^{ss}) \to \overline{G_{\D} (U_{P_K} S)} \cap \overline{G^1}^{ss}$$ sending $(g, z) \to (g, g) \cdot z$ is small and $\overline{G_{\D} (U_{P_K} S)} \cap \overline{G^1}^{ss}$ is a union of strata in $\overline{G^1}^{ss}$. 
\end{thm}

\subsection{} Let $J \subset I$ and $Y$ be a stratum of $L_{J_{\d}} g_0$. By the same argument as above, we can show that $\overline{i_J(G \times_{L_{J_{\d}}} Y/Z(L_J))} \cap \overline{G^1}^{ss}$ is a union of strata of $\overline{G^1}^{ss}$. Since we don't need this result in the rest of the paper, we skip the details.


\section{Character sheaves on $\overline{G^1}^{ss}$}

\subsection{} Fix a prime number $l$ that is invertible in $k$. For any algebraic variety $X$ over $k$, we write $\cd(X)$ for $\cd^b_c(X, \overline{\mathbb Q}_l)$, the bounded derived category of $\overline{\mathbb Q}_l$-constructible sheaves on $X$ (\cite[2.2.18]{BBD}). 

For any subgroup $H$ of $G$ and an $H$-variety $X$, we define the $H$ action on $G \times X$ by $h \cdot (g, x)=(g h \i, h \cdot x)$ and denote by $G \times_H X$ the quotient space. For any perverse sheaf $A$ on $X$ that is equivariant for the $H$ action, we denote by $i^G_H(A)$ the perverse sheaf on $G \times_H X$ such that $p^*(A)[\dim(H)]=\overline{\mathbb Q}_{l, G}[\dim(G)] \boxtimes A$, where $p: G \times X \to G \times_H X$ is the projection map.

\subsection{}\label{xx}
In this subsection, we only assume that $G$ is a connected reductive group. 

Let $\cz=\{g \in Z(G); g g'=g' g \text{ for all } g' \in G^1\}$. For each isolated stratum $S$ of $G^1$ and $n \in \mathbb Z$, let $\cs_n(S)$ be the set of local systems on $S$ that are equivariant for the $\cz^0 \times G$-action defined by $(z, g) \cdot s=g z^n s g \i$ (\cite[5.2]{L2}). Now assume that $\ce$ is an irreducible local system in $\cs_n(S)$.  For $y \in S$, let $H_y$ be the isotropy subgroup of $y$ for this $\cz^0 \times G$-action. Notice that for $(z, g) \in H_y$, $z^n=g \i y g y \i$. By \cite[Lemma 1.1 (2)]{H5}, there are only finitely many possible choices for $z$. In particular, $H_y^0=Z_G(y)^0$. Define a morphism $f: \cz^0 \times G/H^0_y \to S$ by $(z, g) \mapsto (z, g) \cdot y$. Then $\ce$ is a direct summand of $f_! \bbq_{l, \cz \times G/H^0_y}$. Let $C$ be the $L_K$-conjugacy class of $y$, then $f$ factors through \[\xymatrix{\cz^0 \times G/H^0_y \ar[r]^-{f_1} & \cz^0 \times C \ar[r]^-{f_2} & S},\] where $f_1(z, g)=(z^n, g y g \i)$ is a principal $\mu_n \times Z_G(y)/Z_G(y)^0$-covering and $f_2(z, c)=z c$ is a $A$-covering. Here $A=\{z \in \cz^0; z C=C\}$ is a finite group. Therefore $\ce$ is a direct summand of $(f_2)_! (\cf \boxtimes \ce')$, where $\cf$ is an irreducible local system on $\cz^0$ which is a direct summand of $n_! \bbq_{l, \cz^0}$ for the $n$-th isogeny $n: \cz^0 \to \cz^0$ and $\ce'$ is an irreducible local system on $C$ which is a direct summand of $(f_1 \mid_{\{1\} \times G/H^0_y})_! \bbq_{l, G/H^0_y}$. 

Now let $\cz'=\{z \th_{\d}(z) \i; z \in Z(G)^0\}$. Since $Z(G)^0$ is abelian, $\cz'$ is an abelian subgroup of $Z(G)^0$. By \cite[1.2]{L2}, $Z(G)^0=\cz^0 \cz'$. Therefore we have an isomorphism $Z(G)^0 \times_{\cz'} C \cong \cz^0 \times_{\cz^0 \cap \cz'} C$. It is easy to see that $\cz^0 \cap \cz'$ is finite. Since $C \subset G^1$ is stable under the conjugation action of $Z(G)^0$, we have that $\cz' C=C$. Thus we have the following commuting diagram \[\xymatrix{\cz^0 \times \cz' \times C \ar[d]_-a \ar[r]^-b & Z(G)^0 \times C \ar[d]^-c & \\ \cz^0 \times C \ar[r]^-{f_3} & \cz^0 \times_{\cz \cap \cz'} C \ar[r]^-{f_4} & S,}\] where $a, f_3, c$ are projection maps, $b(z, z', c)=(z z', c)$ and $f_4(z, c)=z c$. The square $(a, b, f_3, c)$ is a Cartesian square and $f_2=f_4 \circ f_3$. 

Thus $c^* (f_3)_! (\cf \boxtimes \ce')=b_! a^*(\cf \boxtimes \ce')=b_! (\cf \boxtimes \bbq_{l, \cz'} \boxtimes \ce')$. Any direct summand of $c^* (f_3)_! (\cf \boxtimes \ce')$ is of the form $\cf' \boxtimes \ce'$, where $\cf'$ is an irreducible local system on $Z(G)^0$ which is a direct summand of $n_! \bbq_{l, Z(G)^0}$ for the $n$-th isogeny $n: Z(G)^0 \to Z(G)^0$.
 
As a summary, 

(a) $\ce$ is a direct summand of $(f_4)_! \ce''$. Here $f_4: Z(G)^0 \times_{\cz'} C \to S$, $(z, c) \mapsto zc$ and $\ce''$ is a local system on $Z(G)^0 \times_{\cz'} C$ whose pull back to $Z(G)^0 \times C$ is of the form $\cf' \boxtimes \ce'$, where $\cf'$ is an irreducible local system on $Z(G)^0$ which is a direct summand of $n_! \bbq_{l, Z(G)^0}$ for the $n$-th isogeny $n: Z(G)^0 \to Z(G)^0$.

\begin{lem}\label{1234}
Let $K \subset I$ with $\d(K)=K$, $S \subset L_K g_0$ be an isolated stratum and $S'$ the closure of $S$ in $L_K g_0$. Let $w \in {}^K W \cap W^{\d}$ with $w \i(K) \subset I$. Set $Y=\sqcup_{w \i(K) \subset D \subset I} (S \dot w g_0 \i, \dot w) \cdot h_{D, \d}$ and $Y'=\sqcup_{w \i(K) \subset D \subset I} (S' \dot w g_0 \i, \dot w) \cdot h_{D, \d}$. For $J \subset I$ with $w \i(K) \subset J$, let $\pi_J: S \to (S \dot w g_0 \i, \dot w) \cdot h_{J, \d}$ be the map defined by $s \mapsto (s \dot w g_0 \i, \dot w) \cdot h_{J, \d}$. Let $\ce \in \cs_n(S)$ be an irreducible local system. If $\ce=\pi_J^* \ce'$ for some local system on $Y \cap Z_{J, \d}$, then $IC(Y', \ce) \mid_{Y' \cap Z_{J, \d}}=IC(Y' \cap Z_{J, \d}, \ce')[|I-J|]$. Otherwise, $IC(Y', \ce) \mid_{Y' \cap Z_{J, \d}}=0$. 
\end{lem}

Proof. Let $\tilde Z=\sqcup_{w \i (K) \subset D \subset I} (g_0 Z(L_K) \dot w g_0 \i, \dot w) \cdot h_{D, \d}$ be the closure of $g_0 Z(L_K)$ in $Y'$. Let $p: g_0 Z(L_K) \to \tilde Z \cap Z_{J, \d} \cong g_0 Z(L_K)/\dot w Z(L_J) \dot w \i$, $z \mapsto (z \dot w g_0 \i, \dot w) \cdot h_{J, \d}$ be the projection map . 

We show that 

(a) Let $\cf$ be an irreducible local system on $g_0 Z(L_K)$. If $\cf=p^* \cf'$ for some local system on $\tilde Z \cap Z_{J, \d}$, then $IC(\tilde Z, \cf) \mid_{\tilde Z \cap Z_{J, \d}}=\cf'[|I-J|]$. Otherwise, $IC(\tilde Z, \cf) \mid_{\tilde Z \cap Z_{J, \d}}=0$.

For any $j \notin w \i(K)$, let $\o^\vee_j$ be the fundamental coweight. Then $f_j: \kk^* \to \tilde Z$, $a \mapsto g_0 \dot w \o^\vee_j(a) \dot w \i$ is a cross section to $\tilde Z \cap Z_{I-\{j\}, \d}$ in $\tilde Z$. Using \cite[1.6]{L1}, $IC(\tilde Z, \cf) \mid_{\tilde Z \cap Z_{J, \d}} \neq 0$ if and only if for any $j \notin J$, the monodromy of $\cf$ around the divisor $\tilde Z \cap Z_{I-\{j\}, \d}$ is $0$, i.e., $f_j^* \cf$ is trivial. It is easy to see that $f_j^* \cf$ is trivial for any $j \notin J$ if any only if $\cf=p^* \cf'$. In this case, one can show that $IC(\tilde Z, \cf) \mid_{\tilde Z \cap Z_{J, \d}}=\cf'$. Part (a) is proved. 

Similarly, 

(b) If $\ce=\pi_J^* \ce'$ for some local system on $Y \cap Z_{J, \d}$, then $IC(Y, \ce) \mid_{Y \cap Z_{J, \d}}=\ce'[|I-J|]$. Otherwise, $IC(Y, \ce) \mid_{Y \cap Z_{J, \d}}=0$. 

Let $\cz'=\{z \th_{\d} (z) \i; z \in Z(L_K)\}$. By \ref{xx} (a), $\ce$ is a direct summand of $(f'_4)_! \ce''$. Here $f'_4: g_0 Z(L_K) \times_{\cz'} C \to S$, $(z, c) \mapsto g_0 \i z c=c g_0 \i z$ and $\ce''$ is a local system on $g_0 Z(L_K) \times_{\cz'} C$ whose pull back to $g_0 Z(L_K) \times C$ is of the form $l_{g_0 \i}^* \cf' \boxtimes \ce'$, where $\cf'$ is an irreducible local system on $Z(L_K)$ which is a direct summand of $n_! \bbq_{l, Z(L_K)}$ for the $n$-th isogeny $n: Z(G)^0 \to Z(G)^0$ and $l_{g_0 \i}: g_0 Z(L_K) \to Z(L_K)$, $z \mapsto g_0 \i z$.

Let $C'$ be the closure of $C$ in $L_K g_0$. Then the map $f'_4: g_0 Z(L_K) \times_{\cz'} C \to S$ extends in the natural way to a map $f''_4: \tilde Z \times_{\cz'} C' \to Y'$, $(z, c) \mapsto (c g_0 \i, 1) \cdot z$. This is a surjective map and each fiber is finite. In particular, $f''_4$ is a small map and $$IC(Y', (f'_4)_! \ce'')=(f''_4)_! IC(\tilde Z \times_{\cz'} C', \ce'').$$ 

Consider the following diagram \[\xymatrix{(\tilde Z \cap Z_{J, \d}) \times_{\cz'} C \ar@{^{(}->}[r] \ar[d]_-a & (\tilde Z \cap Z_{J, \d}) \times_{\cz'} C' \ar@{^{(}->}[r] \ar[d]_-b & \tilde Z \times_{\cz'} C' \ar[d]_-{f''_4} \\ Y \cap Z_{J, \d} \ar@{^{(}->}[r] & Y' \cap Z_{J, \d} \ar@{^{(}->}[r] & Y'},\] where $a, b$ are the restriction of $f''_4$ and are small maps. Both squares are Cartesian squares. So $$IC(Y', (f'_4)_! \ce'') \mid_{Y' \cap Z_{J, \d}}=\bigl((f''_4)_! IC(\tz \times_{\cz'} C', \ce'') \bigr) \mid_{Y' \cap Z_{J, \d}}=b_! A,$$ where $A=IC(\tz \times_{\cz'} C', \ce'') \mid_{(\tz \cap Z_{J, \d}) \times_{\cz'} C'}$. 

Notice that the pull back of $A$ to $(\tz \cap Z_{J, \d}) \times C'$ is $IC(\tz, \cf') \mid_{\tz \cap Z_{J, \d}} \boxtimes IC(C', \ce'')$, By (a), the pull back is isomorphic to $$IC((\tz \cap Z_{J, \d}) \times C',  IC(\tz, \cf') \mid_{\tz \cap Z_{J, \d}} \boxtimes \ce'').$$ Here $IC(\tz, \cf') \mid_{\tz \cap Z_{J, \d}} \boxtimes \ce''$ is an irreducible local system on $(\tz \cap Z_{J, \d}) \times C$ or $0$. Hence $A=IC((\tz \cap Z_{J, \d}) \times_{\cz'} C', A \mid_{(\tz \cap Z_{J, \d} \times_{\cz'} C})$ and \begin{align*} IC(Y', (f'_4)_! \ce'') \mid_{Y' \cap Z_{J, \d}} &=b_! A=IC(Y' \cap Z_{J, \d}, a_! (A \mid_{(\tz \cap Z_{J, \d} \times_{\cz'} C})) \\ &=IC(Y' \cap Z_{J, \d}, (b_! A) \mid_{Y \cap Z_{J, \d}}).\end{align*}

Since $IC(Y', \ce)$ is a direct summand of $IC(Y', (f'_4)_! \ce'')$, \begin{align*} IC(Y', \ce) \mid_{Y' \cap Z_{J, \d}} &=IC(Y' \cap Z_{J, \d}, IC(Y', \ce) \mid_{Y \cap Z_{J, \d}}) \\ &=IC(Y' \cap Z_{J, \d}, IC(Y, \ce) \mid_{Y \cap Z_{J, \d}}).\end{align*} Now the lemma follows from (b). \qed

\

From subsection \ref{st} to Lemma \ref{cj}, we only assume that $G$ is a connected reductive group. We first recall some results of character sheaves on disconnected groups. We follow the approach in \cite{L2}. 

\subsection{}\label{st} Let $P$ be a parabolic subgroup of $G$ such that $N_{\tilde G} P \cap G^1 \neq \emptyset$. Let $L$ be a Levi of $P$. Set $L^1=N_{\tilde G} P \cap N_{\tilde G} L \cap G^1$. Consider the diagram \[\xymatrix{L^1 & G \times (N_{\tilde G} P \cap G^1) \ar[l]_-a \ar[r]^-b & G \times_P (N_{\tilde G} P \cap G^1) \ar[r]^-c & G^1,}\] where $a, b$ are projection maps and $c(g, h)=g h g \i$. To any simple perverse sheaf $A$ on $L^1$ which is $L$-equivariant (for the conjugation action) we define $\ind_{L^1}^{G^1} A=c_! A_1$, where $A$ is the perverse sheaf on $G \times_P (N_{\tilde G} P \cap G^1)$ such that $a^* A[\dim(G)-\dim(P)]=b^* A_1$. We call $\ind^{L^1}_{G^1}$ an {\it induction functor}. 

Consider the diagram \[\xymatrix{G^1 & N_{\tilde G} P \cap G^1 \ar[l]_-i \ar[r]^-{\pi} & L^1 \cong (N_{\tilde G} P \cap G^1)/U_P,}\] where $i$ is the inclusion map and $\pi$ is the projection. To any simple perverse sheaf $B$ on $G^1$ which is $G$-equivariant (for the conjugation action), we define $\res^{L^1}_{G^1} B=\pi_! i^* B$. We call $\res^{L^1}_{G^1}$ a {\it restriction functor}. 

\subsection{}\label{5} For $P, L$ and $S$ as in subsection \ref{sm1}, set \begin{gather*} X_{L, S}=G \times_P S' U_P; \\ \tilde Y_{L, S}=G \times_L S^* \cong G \times_P (P_{\D} \cdot S^*). \end{gather*} where $S'$ is the closure of $S$ in $G^1$, $L$ acts diagonally $S^*$ and $P$ acts diagonally on $P_{\D} \cdot S^*$ and $S' U_P$. 

We have the following commuting diagram \[\xymatrix{Y_{L, S} \ar@{^{(}->}[d] & \tilde Y_{L, S} \ar[l]_-{\pi} \ar@{^{(}->}[d] & G \times S^* \ar[l]_-a \ar[r]^-b \ar@{^{(}->}[d] & S \ar@{^{(}->}[d] \\ Y'_{L, S} & X_{L, S} \ar[l]_-{\pi'} & G \times S' U_P \ar[l]_-{a'} \ar[r]^-{b'} & S',}\] where $Y'_{L, S}$ is the closure of $Y_{L, S}$ in $G^1$,  $a, b, a', b'$ are projection maps and $\pi$, $\pi'$ sends $(g, p) \to g p g \i$. 

Let $\ce \in \cs(S)$. Then there is a unique local system $\tilde \ce$ on $\tilde Y_{L, S}$ with $a^* \tilde \ce=b^* \ce$ and the intersection cohomology complex $IC(S', \ce)$, $IC(X_{L, S}, \tilde \ce)$ are related by $(a')^* IC(X_{L, S}, \tilde \ce)=(b')^* IC(S', \ce)$ (see \cite[5.6]{L2}).  Moreover, $IC(Y'_{L, S}, \pi_! \tilde \ce)$ is canonically isomorphic to $\pi'_! IC(X_{L, S}, \tilde \ce)=\ind^{L^1}_{G^1}(IC(S', \ce))[-\dim(X_{L, S})]$ (\cite[Proposition 5.7]{L2}). 

A simple perverse sheaf on $G^1$ is called {\it admissible} if it is a direct summand of the perverse sheaf $IC(Y'_{L, S}, \pi_! \tilde \ce)[\dim(Y'_{L, S})]$ on $G^1$ ($0$ outside $Y'_{L, S}$) for some pair $(L, S)$ as above and a cuspidal local system $\ce \in \cs(S)$ (\cite[6.7]{L2}). 

\begin{lem}\label{cj}
We keep the notations as above. Let $\ce \in \cs(S)$ and $A$ be a direct summand of $IC(Y'_{L, S}, \pi_! \tilde \ce)[\dim(Y'_{L, S})]$. Let $Z$ be a connected subgroup of $Z(G)$. If $A$ is equivariant for the right $Z$-action, then $\ce$ is equivariant for the right $Z$-action on $S$. 
\end{lem}

Proof. Consider the following diagram \[\xymatrix{\tilde Y_{L, S} \ar@{^{(}->}[r] \ar[d]_{\pi} & X_{L, S} \ar[d]^{\pi'} \\ Y_{L, S} \ar@{^{(}->}[r] & Y'_{L, S},}\] where $\pi$ and $\pi'$ are defined in the previous subsection. By \cite[Lemma 5.5]{L2}, this is a Cartesian square. So \begin{align*} \bigl((\pi')^* (\pi')_! IC(X_{L, S}, \tilde \ce) \bigr) \mid_{\tilde Y_{L, S}} &=\pi^* \bigl( (\pi')_! IC(X_{L, S}, \tilde \ce) \mid_{Y_{L, S}} \bigr) \\ &=\pi^* \bigl(IC(Y'_{L, S}, \pi_! \tilde \ce) \mid_{Y_{L, S}} \bigr)=\pi^* \pi_! \tilde \ce. \end{align*}

Consider the following diagram \[\xymatrix{G \times_L (N \times_L S^*) \ar[r]^-b \ar[d]_-a & \tilde Y_{L, S} \ar[d]^-\pi \\ \tilde Y_{L, S} \ar[r]^-\pi & Y_{L, S},}\] where $N=\{n \in N_G L; n S n \i=S\}$ and $G \times_L (N \times_L S^*)$ is the quotient of $G \times (N \times S^*)$ modulo the $L \times L$-action, $(l, l') \cdot (g, n, s)=(g l \i, l n (l') \i, l' s (l') \i)$ and the maps $a$, $b$ are defined by $a(g, n, s)=(g, n s n \i)$ and $b(g, n, s)=(g n, s)$. It is easy to see that this is a Cartesian square. Therefore $\pi^* \pi_! \tilde \ce=b_! a^* \tilde \ce=\tilde \ce^{\oplus |N/L|}$. 

Since $A$ is a direct summand of $(\pi')_! IC(X_{L,S}, \tilde \ce)$, each direct summand of $((\pi')^* A) \mid_{\tilde Y_{L, S}}$ is $\tilde \ce$. In particular, $IC(X_{L, S}, \tilde \ce)$ is an irreducible constitute of ${}^p H^i((\pi')^* A)$ for some $i \in \mathbb Z$. 

Notice that $A$ is equivariant for the right $Z$-action and $\pi'$ is $Z$-equivariant, where the $Z$-action on $X_{L, S}=G \times_P S' U_P$ is defined by $z \cdot (g, s)=(g, s z \i)$. Hence ${}^p H^i((\pi')^* A)$ is also $Z$-equivariant. Therefore $IC(X_{L, S}, \tilde \ce)$ and $IC(X_{L, S}, \tilde \ce) \mid_{G \times_P S U_P}$ are both $Z$-equivariant. By definition, the pull back of $IC(X_{L, S}, \tilde \ce) \mid_{G \times_P S U_P}$ to $G \times S U_P$ is $\bbq_{l, G} \boxtimes \ce \boxtimes \bbq_{l, U_P}$. Therefore $\ce$ is $Z$-equivariant. \qed



\

Now we can prove our main theorem. 

\begin{thm} Let $J, K \subset I$ with $\d(K)=K$. Let $S$ be an isolated stratum of $L_K g_0$, $S'$ its closure in $L_K g_0$ and $\ce \in \cs(S)$. Let $\cw$ be the set of $w \in {}^K W^{J_{\d}} \cap W^{\d}$ with $w \i(K) \subset J_{\d}$ and that $\ce$ is equivariant for the right $\dot w Z(L_J) \dot w \i$-action on $S$. Then $IC(\overline{G^1}^{ss}, \ind^{L_K g_0}_{G^1} IC(S', \ce)) \mid_{Z_{J, \d; 1}}$ is canonically isomorphic to $$\bigoplus_{w \in \cw} i^G_{L_{J_{\d}}} \ind^{L_{w \i(K)} g_0/Z(L_J)}_{L_{J_{\d}} g_0/Z(L_J)} IC({}^{\dot w \i} S'/Z(L_J),  \ce_{J, w}),$$ where $\ce_{J, w}$ is the local system on ${}^{\dot w \i} S/Z(L_J)$ such that $\ce=\Ad(\dot w\i)^* i^* \ce_{J, w}$. Here $i: {}^{\dot w\i} S \to {}^{\dot w \i} S/Z(L_J)$ is the projection map and $\Ad(\dot w \i): S \to {}^{\dot w \i} S$, $s \mapsto \dot w \i s \dot w$. 
\end{thm}

Proof. Consider the following commuting diagram \[\xymatrix{\tilde Y_{L_K, S} \ar@{^{(}->}[r] \ar[d]^{\pi} & X_{L_K, S} \ar@{^{(}->}[r] \ar[d]^{\pi'} & \tilde X_{L_K, S} \ar[d]^{\pi''} \\ Y_{L_K, S} \ar@{^{(}->}[r] & Y'_{L_K, S} \ar@{^{(}->}[r] & \overline{Y'_{L_K, S}} \cap \overline{G^1}^{ss},}\] where $\tilde X_{L_K, S}=G \times_{P_K} (\overline{S U_{P_K}} \cap \overline{G^1}^{ss})$, $\pi''(g, z)=(g, g) \cdot z$ and $\pi, \pi'$ are the restrictions of $\pi''$. Both squares are Cartesian squares. By Theorem \ref{sm2}, $\pi''$ is a small map. Therefore $$IC(\overline{G^1}^{ss}, \ind^{L_K g_0}_{G^1} IC(S', \ce))=IC(\overline{G^1}^{ss}, \pi_! \tilde \ce)[\dim(G)-\dim(L_J)]$$ is canonically isomorphic to $\pi''_! IC(\tx_{L_K, S}, \tilde \ce)[\dim(G)-\dim(L_J)]$. 

Therefore $IC(\overline{G^1}^{ss}, \ind^{L_K g_0}_{G^1} IC(S', \ce))[-\dim(G)+\dim(L_J)] \mid_{Z_{J, \d; 1}}$ is canonically isomorphic to $$\pi''_! IC(\tx_{L_K, S}, \tilde \ce) \mid_{Z_{J, \d; 1}}=(\pi'' \mid_{Z_{J, \d; 1}})_! \bigl(IC(\tx_{L_K, S}, \tilde \ce) \mid_{G \times_{P_K} (\overline{S U_{P_K}} \cap Z_{J, \d; 1})} \bigr).$$

We have shown in Theorem \ref{xxx} that $$\overline{U_{P_K} S} \cap Z_{J, \d; 1}=\sqcup_{w \in {}^K W^{J_{\d}} \cap W^{\d}, w \i(K) \subset J_{\d}} (U_{P_K} S' \dot w g_0 \i, U_{P_K} \dot w) \cdot h_{J, \d}.$$ Similar to the proof of the isomorphism (b) in the proof of Theorem \ref{xxx}, we have that \begin{align*} & (U_{P_K} \cap {}^{\dot w} U) \times (U_{P_K} \cap {}^{\dot w} U^-) \times \sqcup_{J \subset D \subset I} (S' \dot w g_0 \i, U_{P_K} \dot w) \cdot h_{D, \d} \\ & \cong \sqcup_{J \subset D \subset I} (U_{P_K} S' \dot w g_0 \i, U_{P_K} \dot w) \cdot h_{D, \d}.\end{align*} By Lemma \ref{1234}, for $w \in {}^K W^{J_{\d}} \cap W^{\d}$ with $w \i(K) \subset J_{\d}$, the restriction of $IC(\tx_{L_K, S}, \tilde \ce)$ to $G \times_{P_K} (U_{P_K} S' \dot w g_0 \i, U_{P_K} \dot w) \cdot h_{J, \d}$ is $0$ if $w \notin \cw$ and is isomorphic to $\bbq_{(U_{P_K} \cap {}^{\dot w} U) \times (U_{P_K} \cap {}^{\dot w} U^-)} \boxtimes \ce_{J, w}[|I-J|]$ if $w \in \cw$. Now the theorem follows from the isomorphism $G \times_{P_K} (U_{P_K} S' \dot w g_0 \i, U_{P_K} \dot w) \cdot h_{J, \d} \cong G \times_{L_{J_{\d}}} (L_{J_{\d}} \times_{P_{w \i(K)} \cap L_{J_{\d}}} {}^{\dot w \i} S'/Z(L_J))$ in subsection \ref{11x} \& \ref{bbbbb} and subsection \ref{bbbbb} (*). \qed

\subsection{} For any $K \subset J \subset I$ with $\d(K)=K$ and a character sheaf $A$ on $L_K g_0$, we set \[c_J(A)=\begin{cases} A, & \text{ if } A \text{ is equivariant for the right $Z(L_J)$-action on $L_K g_0$}; \\ 0, & \text{ otherwise}. \end{cases}\] If $B$ is a semisimple perverse sheaf on $L_K$ and is a direct sum of some character sheaves $B=\oplus A_i$, then we set $c_J(B)=\oplus c_J(A_i)$. 

By Lemma \ref{cj}, for any $K' \subset K$ with $\d(K')=K'$ and a character sheaf $A$ on $L_{K'} g_0$, if $c_J(A)=0$, then $c_J(\ind^{L_{K'} g_0}_{L_K g_0}(A))=0$. 

Using Macay type formula \cite[Proposition 15.2]{L1} and \cite[Proposition 38.8]{L2}, we can reformulate our main theorem in the following way. 

\begin{thm}\label{ppp}
Let $K \subset I$ with $\d(K)=K$, $S$ be an isolated stratum of $L_K g_0$ and $S'$ be its closure in $L_K g_0$. Let $\ce$ be a cuspidal local system on $S$ and $A=\ind^{L_K g_0}_{G^1} IC(S', \ce[\dim(S)])$. Then for any $J \supset K$, $IC(\overline{G^1}, A) \mid_{Z_{J, \d; 1}}=i^G_{L_{J_{\d}}} (C)[|I-J|]$, where $C$ is a semisimple perverse sheaf on $L_{J_{\d}} g_0/Z(L_J)$ whose pull back to $L_{J_{\d}} g_0$ is $c_J \res^{L_{J_{\d}} g_0}_{G^1} A[-|I-J|]$. 
\end{thm}

By \cite[Section 4]{L1} and \cite[Theorem 30.6]{L2}, any character sheaf on $G$ is a direct summand of $\ind^{L_K g_0}_{G^1} IC(S', \ce[\dim(S)])$ for some pair $(S, \ce)$ as above. We have that 

\begin{cor}
Let $A'$ be a character sheaf on $G^1$ and $J \subset I$. Then $IC(\overline{G^1}, A') \mid_{Z_{J, \d; 1}}$ is of the form $i^G_{L_{\d}}(C)[|I-J|]$ for some semisimple perverse sheaf $C$ on $L_{J_{\d}} g_0/Z(L_J)$.
\end{cor}

Furthermore, we conjecture that the semisimple perverse sheaf $C$ is given by the following explicit formula. 

\begin{conj}\label{conj}
Let $A'$ be a character sheaf on $G^1$. Then for any $J \subset I$, $$IC(\overline{G^1}, A') \mid_{Z_{J, \d; 1}}=i^G_{L_{J_{\d}}} (C)[|I-J|],$$ where $C$ is the semisimple perverse sheaf on $L_{J_{\d}} g_0/Z(L_J)$ whose pull back to $L_{J_{\d}} g_0$ is $c_J \res^{L_{J_{\d}} g_0}_{G^1} (A')[-|I-J|]$. 
\end{conj}

By the above theorem, the conjecture holds for generic character sheaves on $G$. We will show in the end of the paper that this conjecture also holds for Steinberg character sheaf. 

\subsection{} In this subsection, we assume that $G^1=G$. For any $J \subset I$, we have that $\overline{Z_J} \cap \overline{G}^{ss} \cong G \times_{L_J} \overline{G_J}^{ss}$ (see \ref{abs}). Now keep the notation in theorem \ref{ppp}, we can show in the same way as we did for the proof of the main theorem that $$IC(\overline{G}, A) \mid_{\overline{Z_J} \cap \overline{G}^{ss}}=i^G_{L_J} IC(\overline{G_J}^{ss}, C)[|I-J|].$$ Notice that $i^G_{L_J} IC(\overline{G_J}^{ss}, C)$ is canonically isomorphic to $IC(G \times_{L_J} \overline{G_J}^{ss}, i^G_{L_J}(C))$. Thus $IC(\overline{G}, A) \mid_{\overline{Z_J} \cap \overline{G}^{ss}}$ is the intermediate extension of its restriction to $Z_J \cap \overline{G}^{ss}$. Since any character sheaf $A'$ on $G$ is a direct summand of some $A$ considered above, we have that 

(a) For any $J \subset I$, $IC(\overline{G}, A') \mid_{\overline{Z_J} \cap \overline{G}^{ss}}$ is canonically isomorphic to $IC(\overline{Z_J} \cap \overline{G}^{ss}, IC(\overline{G}, A') \mid_{Z_J \cap \overline{G}^{ss}})=i^G_{L_J} IC(\overline{G_J}^{ss}, IC(\overline{G}, A') \mid_{G_J})$.

In particular, for any $K \subset J \subset I$, $IC(\overline{G}, A') \mid_{Z_{K, 1; \d}}$ is canonically isomorphic to $i^G_{L_K} \bigl(IC(\overline{G_J}^{ss}, IC(\overline{G}, A') \mid_{G_J}) \mid_{G_K} \bigr)$. Hence to verify the above conjecture 
for $G^1=G$, it suffices to verify the cases where $J$ is a maximal proper subset of $I$. However, we still don't know how to do it.

Another thing worth mentioning is that the open embedding $G \to \overline{G}^{ss}$ is an affine map. Hence by \cite[Corollary 4.1.12]{BBD}, for any perverse sheaf $A$ on $G$, $IC(\overline{G}^{ss}, A) \mid_{\overline{G}^{ss}-G}[-1]$ is perverse. In other words, $IC(\overline{G}, A) \mid_{\overline{Z_J} \cap \overline{G}^{ss}}[-1]$ is a perverse sheaf for any maximal proper subset $J$ of $I$. We showed above that for any character sheaf $A$, $IC(\overline{G}, A) \mid_{\overline{Z_J} \cap \overline{G}^{ss}}[-|I-J|]$ is perverse for any subset $J$ of $I$. It would be interesting to see if the result holds for arbitrary perverse sheaf on $G$. 

\subsection{} In this and next subsections, we assume that $k$ is an algebraically closure of a finite field $\mathbb F_q$ and that we are given an $\mathbb F_q$-structure on $\tilde G$ with a Frobenius morphism $F: \tilde G \to \tilde G$ such that $G^1$ is defined over $\mathbb F_q$. Then $F$ extends to a Frobenius morphism $F$ on $\overline{G^1}$. 

Let $A$ be a character sheaf on $G^1$ and $\phi: F^* A \to A$ be an isomorphism. Then $\phi$ extends to an isomorphism $F^* IC(\overline{G^1}, A) \to IC(\overline{G^1}, A)$ which we still denote by $\phi$. Then we can define functions $\chi^A_{\phi}: (G^1)^F \to \bbq$ and $\hat \chi^A_{\phi}: (\overline{G^1})^F \to \bbq$ by \begin{gather*} \chi^A_{\phi}(x)=\sum_i (-1)^i Tr(\phi_x^i, H^i(A)_x), \\ \hat \chi^A_{\phi}(x)=\sum_i (-1)^i Tr(\phi_x^i, H^i(IC(\overline{G^1}, A)_x).\end{gather*} The function $\chi^A_{\phi}$ is called the characteristic function of $A$ and is constant on $G^F$-conjugacy classes of $(G^1)^F$ and $\hat \chi^A_{\phi}$ is a natural extension of $\chi^A_{\phi}$. 

Now for any function $f: (\overline{G^1})^F \to \bbq_l$ that is constant on $G^F$-conjugacy classes, we can naturally extend it to a function $\hat f: (\overline{G^1})^F \to \bbq_l$ as follows. 

By \cite[Theorem 25.2]{L1} and \cite[Theorem 21.21]{L2}, the characteristic functions (for various $A$) form a basis of the vector space of functions from $(G^1)^F$ to $\bbq$ that are constant on the $G^F$-conjugacy classes. Hence $f=\sum_A c_A \chi^A_{\phi}$, where $c_A \in \bbq$ is uniquely determined by $f, A$ and $\phi$. Now define \[\hat f=\sum_A c_A \hat \chi^A_{\phi}.\]

We may view the restriction of $\hat f$ to $(\overline{G^1})^F-(G^1)^F$ as the boundary values of $f$.  The most interesting case is when $G^1=G$ and $f$ is an irreducible character of the finite group $G^F$. The study of the boundary values of irreducible characters of $G^F$ is one of the open problems in Springer's talk \cite[Problem 10]{Sp2} at ICM 2006.

\subsection{} The map $(g_1, g_2) \cdot h_{J, \d} \mapsto ({}^{g_2} P_J, {}^{g_1} P_{\d(J)}^-, g_1 U_{P^-_{\d(J)}} g_0 H_{P_J} g_2 \i)$ gives a natural isomorphism of $Z_{J, \d}$ with $\{(P, Q, H_Q g H_P); P \in \cp_J, Q \in \cp_{-w_0(\d(J))}, g \in G^1, {}^g P \cap Q \text{ is a common Levi of } {}^g P \text{ and } Q\}$. Now let $x=(P, Q, \g) \in (\overline{G^1}^{ss})^F$ and $A=IC(\bar G, \ind^{L_K g_0}_{G^1} IC(S', \ce[\dim(S)]))$, where $L_K$ and $S'$ are defined over $F_q$ and $\phi: F^* \ce \to \ce$ is an isomorphism. Then $\phi$ induces a natural isomorphism $F^* A \to A$, which we also denote by $\phi$.  By the previous theorem, we have that \[\tag{*} \hat \chi^A_{\phi}(x)=\sum_{g \in (N_{\tilde G} P \cap \g)^F} \chi^A_{\phi}(g)=\sum_{g \in (N_{\tilde G} Q \cap \g)^F} \chi^A_{\phi}(g).\] 

If the conjecture \ref{conj} is true, then the formula $(*)$ is true for any character sheaf $A$ on $G$ with $\phi: F^* A \cong A$ and $$\hat f(x)=\frac{1}{|(N_{\tilde G} P \cap \g)^F|} \sum_{g \in (N_{\tilde G} P \cap \g)^F} f(g)=\frac{1}{|(N_{\tilde G} Q \cap \g)^F|} \sum_{g \in (N_{\tilde G} Q \cap \g)^F}f(g),$$ for any function $f: (\overline{G^1})^F \to \bbq_l$ that is constant on $G^F$-conjugacy classes and $x=(P, Q, H_Q g H_P) \in (\overline{G^1}^{ss})^F$. 

\

Now we consider a special character sheaf on $G^1$ and its intermediate extension to $\overline{G^1}^{ss}$. 

\subsection{} For any $K \subset J \subset I$ with $\d(K)=K$. The map $L_{J_{\d}} \times (P_K \cap L_{J_{\d}}) g_0/Z(L_J) \to L_{J_{\d}} g_0/Z(L_J)$ defined by $(l, z) \mapsto (l, l) \cdot z$ induces a proper map $$\pi_{J, K, \d}: L_{J_{\d}} \times_{P_K \cap L_{J_{\d}}} (P_K \cap L_{J_{\d}}) g_0/Z(L_J) \to L_{J_{\d}} g_0/Z(L_J).$$ It is known that $\pi_{J, K, \d}$ is a small map. Set $$C_{J, K, \d}=(\pi_{J, K, \d})_! (\overline{\mathbb Q}_{l, L_{J_{\d}} \times_{P_K \cap L_{J_{\d}}} (P_K \cap L_{J_{\d}}) g_0/Z(L_J)}[\dim(G_{J_{\d}})]).$$

Moreover, we may identify $(L_{J_{\d}}, 1) \cdot h_{J, \d}$ with $L_{J_{\d}} g_0/Z(L_J)$ and $(P_K \cap L_{J_{\d}}, 1) \cdot h_{J, \d}$ with $(P_K \cap L_{J_{\d}}) g_0/Z(L_J)$ in the natural way. Under this identification, $C_{J, K, \d}$ is a perverse sheaf on $(L_{J_{\d}}, 1) \cdot h_{J, \d}$. 

Define $\pi'_{J, K, \d}: G \times_{P_K \cap L_{J_{\d}}} (P_K \cap L_{J_{\d}}, 1) \cdot h_{J, \d} \to Z_{J, 1; \d}$ by $(g, z) \mapsto (g, g) \cdot z$. Notice that $$G \times_{L_{J_{\d}}} \bigl(L_{J_{\d}} \times_{P_K \cap L_{J_{\d}}} (P_K \cap L_{J_{\d}}, 1) \cdot h_{J, \d} \bigr) \cong G \times_{P_K \cap L_{J_{\d}}} (P_K \cap L_{J_{\d}}, 1) \cdot h_{J, \d}.$$ 

Then \[i^G_{L_{J_{\d}}}(C_{J, K, \d})=(\pi'_{J, K, \d})_! (\overline{\mathbb Q}_{l, G \times_{P_K \cap L_{J_{\d}}} (P_K \cap L_{J_{\d}}, 1) \cdot h_{J, \d}}[\dim(Z_{J, 1; \d})])\] is a perverse sheaf on $G \times_{L_{J_{\d}}} (L_{J_{\d}}, 1) \cdot h_{J, \d} \cong Z_{J, 1; \d}$. 

\subsection{} Let $J, K \subset I$ with $\d(K)=K$ and $w \in {}^K W^{J_{\d}} \cap W^{\d}$. Let $\e(w)=\min(w W_J)$. Set $K_1=\max\{K' \subset K; \d(K')=K', \e(w) \i(K') \subset J\}$. By Lemma \ref{7}, $K_1=K \cap w J_{\d}$. By subsection \ref{p} (7), \begin{align*} (P_K)_{\D} (B \dot w, B \dot w) \cdot h_{J, \d} &=[J, \d(\e(w)), \e(w)]_{K, \d} \\ &=(P_K)_{\D} (L_{K_1} \dot{\d(\e(w))}, \dot{\e(w)}) \cdot h_{J, \d} \\ &=(P_K)_{\D} (L_{K_1} \dot w, \dot w) \cdot h_{J, \d} \\ &=(P_K)_{\D} (\dot w L_{w \i(K) \cap J_{\d}}, \dot w) \cdot h_{J, \d}. \end{align*}

The map $$f: G \times_{P_K} (P_K)_{\D} (B \dot w, B \dot w) \cdot h_{J, \d}  \to G \times_{{}^{\dot w \i} P_K} ({}^{\dot w \i} P_K)_{\D} (L_{w \i(K) \cap J_{\d}}, 1) \cdot h_J$$ defined by $(g, z) \mapsto (g \dot w, (\dot w \i, \dot w \i) z)$ is an isomorphism. Moreover, $\pi_{J, K, w}=\pi'_{J, K, w} \circ f$, where \begin{gather*} \pi_{J, K, w, \d}: G \times_{P_K} (P_K)_{\D} (B \dot w, B \dot w) \cdot h_{J, \d}  \to Z_{J, 1; \d}, \\ \pi'_{J, K, w, \d}:  G \times_{{}^{\dot w \i} P_K} ({}^{\dot w \i} P_K)_{\D} (L_{w \i(K) \cap J_{\d}}, 1) \cdot h_J \to Z_{J, 1; \d}, \end{gather*} are induced from the map $G \times Z_{J, 1; \d} \to Z_{J, 1; \d}$ defined by $(g, z) \mapsto (g, g) \cdot z$.

Notice that $$(P_{w \i(K) \cap J_{\d}} \cap L_{J_{\d}})_{\D} (L_{w \i(K) \cap J_{\d}}, 1) \cdot h_{J, \d}=(P_{w \i(K) \cap J_{\d}} \cap L_{J_{\d}}, 1) \cdot h_{J, \d}.$$ By subsection \ref{11x} (a), $\pi_{J, K, w, \d}$ is a small map and \begin{align*} & (\pi_{J, K, w, \d})_! (\overline{\mathbb Q}_{l, G \times_{P_K} (P_K)_{\D} (B \dot w, B \dot w) \cdot h_{J, \d}} [\dim(Z_{J, 1; \d})]) \\ & =(\pi'_{J, K, w, \d})_! (\overline{\mathbb Q}_{l, G \times_{{}^{\dot w \i} P_K} ({}^{\dot w \i} P_K)_{\D} (L_{w \i(K) \cap J_{\d}}, 1) \cdot h_{J, \d}} [\dim(Z_{J, 1; \d})])\\ &=(\pi'_{J, w \i(K) \cap J_{\d}, \d})_! (\overline{\mathbb Q}_{l, G \times_{P_{w \i(K) \cap J_{\d}} \cap L_{J_{\d}}} (P_{w \i(K) \cap J_{\d}} \cap L_{J_{\d}}, 1) \cdot h_{J, \d}}[\dim(Z_{J, 1; \d})]) \\ &=i^G_{L_{J_{\d}}}(C_{J, w \i(K) \cap J_{\d}, \d}). \end{align*}

\subsection{} For any $J \subset I$ with $\d(J)=J$, there is a unique simple perverse sheaf $St_{J, \d}$ on $(L_{J_{\d}}, 1) \cdot h_{J, \d} \cong L_{J_{\d}} g_0/Z(L_J)$ such that $St_{J, \d}$ is a direct summand of $C_{J, \emptyset, \d}$ and $St_{J, \d}$ is not a direct summand of $C_{J, K, \d}$ for any $\emptyset \neq K \subset J$ with $\d(K)=K$. In fact, $$St_{J, \d} \bigoplus \bigoplus\limits_{K \subset J, \d(K)=K, 2 \nmid |K|} C_{J, K, \d}=\bigoplus\limits_{K \subset J, \d(K)=K, 2 \mid |K|} C_{J, K, \d}.$$

It is known that for any $g \in G_J^1$, $\mathcal H^i_g(St_{J, \d}) \neq 0$ for some $i \in \mathbb Z$ if and only if the stabilizer of $g$ in $G$ is reductive (i.e., $g$ is quasi-semisimple). In this case, $\sum_{i \in \mathbb Z} \dim(\mathcal H^i_g(St_{J, \d}))=1$.  See \cite[12.6]{L2}.

Let $K \subset I$ with $\d(K)=K$. By Theorem \ref{smooth}, $\overline{P_K} \cap \overline{G^1}^{ss}$ is smooth. By Theorem \ref{bp} and the previous subsection, $\pi_K: G \times_{P_K} (\overline{P_K} \cap \overline{G^1}^{\ss}) \to \overline{G^1}^{ss}$ defined by $(g, z) \mapsto (g, g) \cdot z$ is a small map. Hence $(\pi_K)_! (\overline{\mathbb Q}_{l, G \times_{P_K} (\overline{P_K} \cap \overline{G^1}^{\ss})} [\dim(G)])$ is a perverse sheaf on $\overline{G^1}^{ss}$ whose restriction to $G^1$ is $C_{K, I, \d}$. 

Let $S'$ be the unique simple perverse sheaf on $\overline{G^1}^{ss}$ such that $S' \mid_{G^1}=St_{I, \d}$. Then $S' \bigoplus \bigoplus\limits_{K \subset I, \d(K)=K, 2 \nmid |K|} (\pi_K)_! (\overline{\mathbb Q}_{l, G \times_{P_K} (\overline{P_K} \cap \overline{G^1}^{\ss})} [\dim(G)]) \linebreak=\bigoplus\limits_{K \subset I, \d(K)=K, 2 \mid |K|} (\pi_K)_! (\overline{\mathbb Q}_{l, G \times_{P_K} (\overline{P_K} \cap \overline{G^1}^{\ss})} [\dim(G)])$.

\

Now we calculate the restriction of $S'$ to $Z_{J, 1; \d}$. 

\begin{prop} For $J \subset I$, $S' \mid_{Z_{J, 1; \d}}=i^G_{L_{J_{\d}}}(St_{J, {\d}}[|I-J|])$.
\end{prop}

Proof. For $J, K \subset I$ and $w \in {}^K W^{J_{\d}} \cap W^{\d}$, set $I(J, K, w, \d)=w \i K \cap J_{\d}$. 

We have that \begin{align*} & \bigl((\pi_K)_! (\overline{\mathbb Q}_{l, G \times_{P_K} (\overline{P_K} \cap \overline{G^1}^{\ss})} [\dim(G)])\bigr) \mid_{Z_{J, 1; \d}} \\ &=(\pi_K \mid_{G \times_{P_K} (\overline{P_K} \cap Z_{J, 1; \d})})_! (\overline{\mathbb Q}_{l, G \times_{P_K} (\overline{P_K} \cap Z_{J, 1; \d})} [\dim(G)]) \\ &=(\pi_K \mid_{G \times_{P_K} (\overline{P_K} \cap Z_{J, 1; \d})})_! \overline{\mathbb Q}_{l, \sqcup_{w \in {}^K W^{J_{\d}} \cap W^{\d}} G \times_{P_K} (P_K)_{\D} (B \dot w, B \dot w) \cdot h_{J, \d}} [\dim(G)] \\ &=\bigoplus_{w \in {}^K W^{J_{\d}} \cap W^{\d}}  (\pi_{J, K, w, \d})_! (\overline{\mathbb Q}_{l, G \times_{P_K} (P_K)_{\D} (B \dot w, B \dot w) \cdot h_{J, \d}} [\dim(G)]) \\ &=\bigoplus_{w \in {}^K W^{J_{\d}} \cap W^{\d}} i^G_{L_{J_{\d}}} (C_{J, I(J, K, w, \d), \d})[|I-J|]. \end{align*}

Moreover, \begin{align*} & \bigoplus_{K \subset I, \d(K)=K, 2 \mid |K|} \bigoplus_{w \in {}^K W^{J_{\d}} \cap W^{\d}} i^G_{L_{J_{\d}}} (C_{J, I(J, K, w, \d), \d})[|I-J|] \\ &=\bigoplus_{w \in W^{J_{\d}} \cap W^{\d}} \bigoplus_{K \subset I, \d(K)=K, 2 \mid |K|, w \in {}^K W} i^G_{L_{J_{\d}}} (C_{J, I(J, K, w, \d), \d})[|I-J|] \\ &=\bigoplus_{w \in W^{J_{\d}} \cap W^{\d}} \bigoplus_{K \subset I(J, I, w, \d)} i^G_{L_{J_{\d}}} (C_{J, K, \d})[|I-J|] \bigoplus_{K' \subset I, \d(K')=K', w \in {}^{K'} W, I(J, K', w, \d)=K, 2 \mid |K'|} 1.\end{align*} 

Similarly, \begin{align*} & \bigoplus_{K \subset I, \d(K)=K, 2 \nmid |K|} \bigoplus_{w \in {}^K W^{J_{\d}} \cap W^{\d}} i^G_{L_{J_{\d}}} (C_{J, I(J, K, w, \d), \d})[|I-J|] \\ &=\bigoplus_{w \in W^{J_{\d}} \cap W^{\d}} \bigoplus_{K \subset I(J, I, w, \d)} i^G_{L_{J_{\d}}} (C_{J, K, \d})[|I-J|] \bigoplus_{K' \subset I, \d(K')=K', w \in {}^{K'} W, I(J, K', w, \d)=K, 2 \nmid |K'|} 1.\end{align*} 


Fix $w \in W^{J_{\d}} \cap W^{\d}$. Let $J'=\max\{K \subset I; w \in {}^K W\}$. Then $\d(J')=J'$ and $w I(J, I, w, \d) \subset J'$. It is easy to see that for any $K \subset I(J, I, w, \d)$ and $K' \subset I$ with $\d(K')=K'$, the following conditions are equivalent:

(1) $w \in {}^{K'} W$ and $I(J, K', w, \d)=K$;

(2) $K=\d(K)$ and $w K \subset K' \subset w K \sqcup (J'-I(J, I, w, \d))$.

Notice that $J'-I(J, I, w, \d)$ is $\d$-stable. Therefore, for any $K \subset I(J, I, w, \d)$, \[\sum_{K' \subset I, w \in {}^{K'} W, I(J, K', w, \d)=K} (-1)^{|K'|}=\begin{cases} (-1)^{|K|}, & \text{ if } \d(K)=K \text{ and } J'-I(J, I, w, \d)=\emptyset; \\ 0, & \text{ otherwise}. \end{cases}\]

If $J'=I(J, I, w, \d)$, then there is no $i \in I$ such that $w w_0^{J_{\d}} \in {}^{\{i\}} W$. Hence $w w_0^{J_{\d}}=w_0$ and $w=w_0 w_0^{J_{\d}}$. In this case, $I(J, I, w, \d)=J_{\d}$. Therefore for any $w \in W^{J_{\d}} \cap W^{\d}$ with $w \neq w_0 w_0^{J_{\d}}$ \begin{align*} & \bigoplus_{K \subset I(J, I, w, \d)} i^G_{L_{J_{\d}}} (C_{J, K, \d})[|I-J|] \bigoplus_{K' \subset I, \d(K')=K', w \in {}^{K'} W, I(J, K', w, \d)=K, 2 \mid |K'|} 1 \\ &=\bigoplus_{K \subset I(J, I, w, \d)} i^G_{L_{J_{\d}}} (C_{J, K, \d})[|I-J|] \bigoplus_{K' \subset I, \d(K')=K', w \in {}^{K'} W, I(J, K', w, \d)=K, 2 \nmid |K'|} 1. \end{align*}

Now \begin{align*} S' \mid_{Z_{J, 1; \d}} \bigoplus & \bigoplus_{K \subset J_{\d}, \d(K)=K, 2 \nmid |K| } i^G_{L_{J_{\d}}} (C_{J, K, \d})[|I-J|] \\ &=\bigoplus_{K \subset J_{\d}, \d(K)=K, 2 \mid |K|} i^G_{L_{J_{\d}}} (C_{J, K, \d})[|I-J|]. \end{align*} Hence $S' \mid_{Z_{J, 1; \d}}=i^G_{L_{J_{\d}}}(St_{J, \d}[|I-J|])$.  \qed

\subsection{} Let $\tilde S$ be the simple perverse sheaf on $\overline{G^1}$ such that $\tilde S \mid_{G^1}=St_{I, \d}$. Then $\tilde S \mid_{\overline{G^1}^{ss}}=S'$. By the previous Proposition and \cite{HS}, the following conditions on $z \in \overline{G^1}$ are equivalent:

(1) The stabilizer of $z$ in $G$ is reductive;

(2) $z \in \overline{G^1}^{ss}$ and $\mathcal H^i_z(\tilde S) \neq 0$ for some $i \in \mathbb Z$;

(3) $z \in \overline{G^1}^{ss}$ and $\sum_{i \in \mathbb Z} \dim(\mathcal H^i_z(\tilde S))=1$.

This verifies Lusztig's conjecture in \cite[12.6]{L3} inside $\overline{G^1}^{ss}$. More precisely, by what we have shown above, Lusztig's conjecture is now reduced to the following one:

\begin{conj} The intermediate extension of $S'$ to $\overline{G^1}$ is the extension by $0$ outside $\overline{G^1}^{ss}$. 
\end{conj}

\section*{Acknowledgement} 
We thank George Lusztig for his continuous encouragement and many helpful discussions on characters sheaves. We thank T. A. Springer for valuable suggestions on a preliminary version of the paper, especially for the study of the restriction of $IC(\overline{G}, A)$ to the closure of $Z_{J, 1; id}$ in $\overline{G}^{ss}$ instead of $Z_{J, 1; id}$. 

\bibliographystyle{amsalpha}

\end{document}